\documentclass[11pt,leqno]{article}
\usepackage{amssymb}
\usepackage[mathscr]{eucal}
\usepackage{amsmath,amssymb,latexsym,theorem,bbm}
\usepackage{color}
\usepackage{appendix}

\newcommand{\NN}{\mathbb{N}}

\newcommand{\RR}{\mathbb{R}}

\newcommand{\ZZ}{\mathbb{Z}}

\newcommand{\bA}{{\boldsymbol{A}}}

\newcommand{\bD}{{\boldsymbol{D}}}

\newcommand{\bE}{{\boldsymbol{E}}}

\newcommand{\bI}{{\boldsymbol{I}}}
\newcommand{\bJ}{{\boldsymbol{J}}}

\newcommand{\bM}{{\boldsymbol{M}}}

\newcommand{\bx}{{\boldsymbol{x}}}

\newcommand{\bzero}{{\boldsymbol{0}}}

\newcommand{\cB}{{\mathcal B}}

\newcommand{\cF}{{\mathcal F}}

\newcommand{\cL}{{\mathcal L}}

\newcommand{\cN}{{\mathcal N}}

\newcommand{\dd}{\mathrm{d}}
\newcommand{\ee}{\mathrm{e}}

\newcommand{\CLSED}{\mathrm{CLSE}}

\DeclareMathOperator*{\argmin}{arg\,min}

\newcommand{\EE}{\operatorname{\mathbb{E}}}
\newcommand{\PP}{\operatorname{\mathbb{P}}}

\newcommand{\ha}{\widehat{a}}
\newcommand{\hb}{\widehat{b}}
\newcommand{\hc}{\widehat{c}}
\newcommand{\hd}{\widehat{d}}

\newcommand{\halpha}{\widehat{\alpha}}
\newcommand{\hbeta}{\widehat{\beta}}
\newcommand{\hgamma}{\widehat{\gamma}}
\newcommand{\hdelta}{\widehat{\delta}}

\newcommand{\tW}{\widetilde{W}}

\newcommand{\vare}{\varepsilon}

\renewcommand{\mid}{\,|\,}

\renewcommand{\leq}{\leqslant}
\renewcommand{\geq}{\geqslant}

\newcommand{\stoch}{\stackrel{\PP}{\longrightarrow}}
\newcommand{\distr}{\stackrel{\cL}{\longrightarrow}}

\newcommand{\as}{\stackrel{{\mathrm{a.s.}}}{\longrightarrow}}

\newcommand{\bbone}{\mathbbm{1}}

\newcommand{\proofend}{\hfill\mbox{$\Box$}}

\numberwithin{equation}{section}

\theoremstyle{change} \theorembodyfont{\em}
\newtheorem{Lem}{Lemma.}[section]
\newtheorem{Thm}[Lem]{Theorem.}
\newtheorem{Pro}[Lem]{Proposition.}

\newtheorem{Def}[Lem]{Definition.}

\theorembodyfont{\rm}
\newtheorem{Rem}[Lem]{Remark.}

\begin{document}

\begin{center}
 {\bfseries\Large
  Parameter estimation for the subcritical Heston model\\[2mm]
 based on discrete time observations} \\[5mm]
 {\sc\large
  M\'aty\'as $\text{Barczy}^{*,\diamond}$,
  \ Gyula $\text{Pap}^{**}$
  \ and
   Tam\'as $\text{T. Szab\'o}^{**}$}
\end{center}

\vskip0.2cm

\noindent
 * Faculty of Informatics, University of Debrecen,
   Pf.~12, H--4010 Debrecen, Hungary.

\noindent
 ** Bolyai Institute, University of Szeged,
      Aradi v\'ertan\'uk tere 1, H--6720 Szeged, Hungary.

\noindent e--mails: barczy.matyas@inf.unideb.hu (M. Barczy),
                    papgy@math.u-szeged.hu (G. Pap),
                    tszabo@math.u-szeged.hu (T. T. Szab\'o).

\noindent $\diamond$ Corresponding author.

\renewcommand{\thefootnote}{}
\footnote{\textit{2010 Mathematics Subject Classifications\/}:
    91G70, 60H10, 62F12, 60F05.}
\footnote{\textit{Key words and phrases\/}:
 Heston model, conditional least squares estimation.}
\vspace*{0.2cm}
\footnote{This research was supported by the European Union and the State of
 Hungary, co-financed by the European Social Fund in the framework of
 T\'AMOP-4.2.4.A/ 2-11/1-2012-0001  ‘National Excellence Program’.}

\vspace*{-10mm}

%\centerline{\sl \today}

\begin{abstract}
We study asymptotic properties of some (essentially conditional least squares)
 parameter estimators for the subcritical Heston model based on discrete time
 observations derived from conditional
 least squares estimators of some modified parameters.
\end{abstract}

\section{Introduction}

The Heston model has been extensively used in financial mathematics since one
 can well-fit them to real financial data set, and they are well-tractable
 from the point of view of computability as well.
Hence parameter estimation for the Heston model is an important task.

In this paper we study the Heston model
 \begin{align}\label{Heston_SDE}
  \begin{cases}
   \dd Y_t = (a - b Y_t) \, \dd t + \sigma_1 \sqrt{Y_t} \, \dd W_t , \\
   \dd X_t = (\alpha - \beta Y_t) \, \dd t
             + \sigma_2  \sqrt{Y_t}
               \bigl(\varrho \, \dd W_t
                     + \sqrt{1 - \varrho^2} \, \dd B_t\bigr) ,
  \end{cases} \qquad t \geq 0 ,
 \end{align}
 where \ $a > 0$, \ $b, \alpha, \beta \in \RR$, \ $\sigma_1 > 0$,
 \ $\sigma_2 > 0$, \ $\varrho \in (-1, 1)$, \ and \ $(W_t, B_t)_{t\geq0}$ \ is a
 2-dimensional standard Wiener process, see Heston \cite{Hes}.
We investigate only the so-called subcritical case, i.e., when \ $b > 0$,
 \ see Definition \ref{Def_criticality}, and we introduce some parameter
 estimator of \ $(a, b, \alpha, \beta)$ \ based on discrete time observations
 and derived from conditional least squares estimators (CLSEs) of some
 modified parameters starting the process \ $(Y,X)$ \ from some known non-random initial
 value \ $(y_0,x_0)\in(0,\infty)\times\RR$.
\ We do not estimate the parameters \ $\sigma_1$, $\sigma_2$ \ and \ $\varrho$,
 \ since these parameters could---in principle, at least---be
 determined (rather than estimated) using an arbitrarily short continuous time
 observation \ $(X_t)_{t\in[0,T]}$ \ of \ $X$, \ where $T>0$, \ see, e.g., Barczy and Pap
 \cite[Remark 2.6]{BarPap2}.
In Overbeck and Ryd\'en \cite[Theorems 3.2 and 3.3]{OveRyd}
 one can find a strongly consistent and asymptotically normal estimator of
 \ $\sigma_1$ \ based on discrete time observations for the process \ $Y$,
 \ and for another estimator of \ $\sigma_1$, \ see Dokuchaev \cite{Dok}.
\ Eventually, it turns out that for the calculation of the estimator of
 \ $(a, b, \alpha, \beta)$, \ one does not need to know the values of the
 parameters \ $\sigma_1, \sigma_2$ \ and \ $\varrho$.
\ For interpretations of \ $Y$ \ and \ $X$ \ in financial mathematics, see, e.g.,
 Hurn et al. \cite[Section 4]{HurLinMcC}.

CLS estimation has been considered for the Cox-Ingersoll-Ross (CIR) model,
 which satisfies the first equation of \eqref{Heston_SDE}.
For the CIR model, Overbeck and Ryd\'en \cite{OveRyd} derived the CLSEs and
 gave their asymptotic properties, however, they did not investigate the
 conditions of their existence.
Specifically, Theorems 3.1 and 3.3 in Overbeck and Ryd\'en \cite{OveRyd}
 correspond to our Theorem \ref{Thm_CLSE_a_b_alpha_beta}, but they estimate
 the volatility coefficient \ $\sigma_1$ \ as well, which we assume to be
 known.
Li and Ma \cite{LiMa} extended the investigation to so-called stable CIR
 processes driven by an $\alpha$-stable process instead of a Brownian motion.
For a more complete overview of parameter estimation for the Heston model see,
 e.g., the introduction in Barczy and Pap \cite{BarPap2}.

It would be possible to calculate the discretized version of the maximum
 likelihood estimators derived in Barczy and Pap
 \cite{BarPap2} using the same procedure as in Ben Alaya and Kebaier
 \cite[Section 4]{BenKeb2} valid for discrete time observations of high
 frequency.
However, this would be basically different from the present line of
 investigation, therefore we will not discuss it further.

The organization of the paper is the following.
In Section 2 we recall some important results about the existence of a unique
 strong solution to \eqref{Heston_SDE}, and study its asymptotic properties.
In the subcritical case, i.e., when \ $b>0$, \ we invoke a result due to
 Cox et al. \cite{CoxIngRos} on the unique existence of a stationary
 distribution, and we slightly improve a result due to Li and Ma \cite{LiMa}
 and Jin et al. \cite[Corollary 2.7]{JinManRudTra} and \cite[Corollaries 5.9 and 6.4]{JinRudTra}
 on the ergodicity of the CIR process \ $(Y_t)_{t\geq 0}$,
 \ see Theorem \ref{Ergodicity}.
We also recall some convergence results for square-integrable martingales.
In Section 3 we introduce the CLSE of a transformed parameter vector based on
 discrete time observations, and derive the asymptotic properties of the
 estimates -- namely, strong consistency and asymptotic normality, see
 Theorem \ref{c_d_normal}.
Thereafter, we apply these results together with the so-called delta method to
 obtain the same asymptotic properties of the estimators for the original
 parameters, see Theorem \ref{Thm_CLSE_a_b_alpha_beta}.
The point of the parameter transformation is to reduce the minimization in the
 CLS method to a linear problem,  because our objective function depends on
 the original parameters through complicated functions.
The covariance matrices of the limit normal distributions in Theorems \ref{c_d_normal}
 and \ref{Thm_CLSE_a_b_alpha_beta} depend on the unknown parameters \ $a$, $b$ \ and
 \ $\beta$, \ as well (but somewhat surprisingly not on \ $\alpha$).
They also depend on the
 volatility parameters $\sigma_1$, $\sigma_2$ and $\rho$, but, again, we will assume these to be known.
Since the considered estimators of \ $a$, $b$ \ and \ $\beta$ are proved to be strongly
 consistent, using random normalization, one may derive counterparts of Theorems
 \ref{c_d_normal} and \ref{Thm_CLSE_a_b_alpha_beta} in a way that the limit distributions
 are four-dimensional standard normal distributions (having the identity matrix
 \ $\bI_4$ \ as covariance matrices).

\section{Preliminaries}
\label{Prel}

Let \ $\NN$, \ $\ZZ_+$, \ $\RR$, \ $\RR_+$, \ $\RR_{++}$, \ and
 \ $\RR_{--}$ \ denote the sets of positive integers, non-negative integers,
 real numbers, non-negative real numbers, positive real numbers, and negative real numbers, respectively.
For \ $x , y \in \RR$, \ we will use the notation \ $x \land y := \min(x, y).$
\ By \ $\|\bx\|$ \ and \ $\|\bA\|$, \ we denote the Euclidean norm of a vector
 \ $\bx \in \RR^d$ \ and the induced matrix norm of a matrix
 \ $\bA \in \RR^{d \times d}$, \ respectively.
By \ $\bI_d \in \RR^{d \times d}$, \ we denote the \ $d\times d$ unit matrix.
The Borel \ $\sigma$-algebra on \ $\RR$ \ is denoted by \ $\cB(\RR)$.
\ Let \ $\bigl(\Omega, \cF, \PP\bigr)$ \ be a probability space
 equipped with the augmented filtration \ $(\cF_t)_{t\in\RR_+}$ \ corresponding to
 \ $(W_t,B_t)_{t\in\RR_+}$ \ and a given initial value \ $(\eta_0,\zeta_0)$ \ being independent
 of \ $(W_t,B_t)_{t\in\RR_+}$ \ such that \ $\PP(\eta_0\in\RR_+)=1$, \
 constructed as in Karatzas and Shreve \cite[Section 5.2]{KarShr}.
Note that \ $(\cF_t)_{t\in\RR_+}$ \ satisfies the usual conditions, i.e., the
 filtration \ $(\cF_t)_{t\in\RR_+}$ \ is right-continuous and \ $\cF_0$
 \ contains all the $\PP$-null sets in \ $\cF$.

The next proposition is about the existence and uniqueness of a strong
 solution of the SDE \eqref{Heston_SDE}, see, e.g., Barczy and Pap
 \cite[Proposition 2.1]{BarPap2}.

\begin{Pro}\label{Pro_Heston}
Let \ $(\eta_0, \zeta_0)$ \ be a random vector independent of
 \ $(W_t, B_t)_{t\in\RR_+}$ \ satisfying \ $\PP(\eta_0 \in \RR_+) = 1$.
\ Then for all \ $a \in \RR_{++}$, \ $b, \alpha, \beta \in \RR$,
 \ $\sigma_1, \sigma_2 \in \RR_{++}$, \ and \ $\varrho \in (-1, 1)$,
 \ there is a (pathwise) unique strong solution \ $(Y_t, X_t)_{t\in\RR_+}$ \ of
 the SDE \eqref{Heston_SDE} such that
 \ $\PP((Y_0, X_0) = (\eta_0, \zeta_0)) = 1$ \ and
 \ $\PP(\text{$Y_t \in \RR_+$ \ for all \ $t \in \RR_+$}) = 1$.
\ Further, for all \ $s, t \in \RR_+$ \ with \ $s \leq t$,
 \begin{align}\label{Solutions}
  \begin{cases}
   Y_t = \ee^{-b(t-s)} Y_s
         + a \int_s^t \ee^{-b(t-u)} \, \dd u
         + \sigma_1 \int_s^t \ee^{-b(t-u)} \sqrt{Y_u} \, \dd W_u, \\
   X_t = X_s + \int_s^t (\alpha - \beta Y_u) \, \dd u
         + \sigma_2
           \int_s^t
            \sqrt{Y_u} \, \dd (\varrho W_u + \sqrt{1 - \varrho^2} B_u) .
  \end{cases}
 \end{align}
\end{Pro}

Next we present a result about the first moment of \ $(Y_t, X_t)_{t\in\RR_+}$.
\ For a proof, see, e.g., Barczy and Pap \cite[Proposition 2.2]{BarPap2} together
  with \eqref{Solutions} and Proposition 3.2.10 in Karatzas and Shreve \cite{KarShr}.

\begin{Pro}\label{Pro_moments}
Let \ $(Y_t, X_t)_{t\in\RR_+}$ \ be the unique strong solution of the SDE
 \eqref{Heston_SDE} satisfying \ $\PP(Y_0 \in \RR_+) = 1$ \ and
 \ $\EE(Y_0) < \infty$, \ $\EE(|X_0|) < \infty$.
\ Then for all \ $s, t \in \RR_+$ \ with \ $s \leq t$, \ we have
 \begin{align}\label{cond_exp_discrete_Y}
  &\EE(Y_t \mid \cF_s)
     = \ee^{-b(t-s)} Y_s + a \int_s^t \ee^{-b(t-u)} \, \dd u,\\  \label{cond_exp_discrete_X}
  &\EE(X_t \mid \cF_s)
    = X_s + \int_s^t (\alpha - \beta \EE(Y_u \mid \cF_s)) \, \dd u \\\nonumber
  &\phantom{\EE(X_t \mid \cF_s)\,}
     = X_s + \alpha (t - s)
       - \beta Y_s \int_s^t \ee^{-b(u-s)}\,\dd u
       - a \beta \int_s^t \left(\int_s^u \ee^{-b(u-v)} \, \dd v\right) \dd u ,
 \end{align}
 and hence
 \begin{align*}
  \begin{bmatrix}
   \EE(Y_t) \\
   \EE(X_t) \\
  \end{bmatrix}
  = \begin{bmatrix}
     \ee^{-bt} & 0 \\
     - \beta \int_0^t \ee^{-bu} \, \dd u & 1 \\
    \end{bmatrix}
    \begin{bmatrix}
     \EE(Y_0) \\
     \EE(X_0) \\
    \end{bmatrix}
    + \begin{bmatrix}
       \int_0^t \ee^{-bu} \, \dd u & 0 \\
       - \beta \int_0^t \left(\int_0^u \ee^{-bv} \, \dd v \right) \dd u & t \\
      \end{bmatrix}
      \begin{bmatrix}
       a \\
       \alpha \\
      \end{bmatrix} .
 \end{align*}
Consequently, if \ $b \in \RR_{++}$, \ then
 \[
   \lim_{t\to\infty} \EE(Y_t) = \frac{a}{b} , \qquad
   \lim_{t\to\infty} t^{-1} \EE(X_t) = \alpha - \frac{\beta a}{b} ,
 \]
 if \ $b = 0$, \ then
 \[
   \lim_{t\to\infty} t^{-1} \EE(Y_t) = a , \qquad
   \lim_{t\to\infty} t^{-2} \EE(X_t) = - \frac{1}{2} \beta a ,
 \]
 if \ $b \in \RR_{--}$, \ then
 \[
   \lim_{t\to\infty} \ee^{bt} \EE(Y_t) = \EE(Y_0) - \frac{a}{b} , \qquad
   \lim_{t\to\infty} \ee^{bt} \EE(X_t)
   = \frac{\beta}{b} \EE(Y_0) - \frac{\beta a}{b^2} .
 \]
\end{Pro}

Based on the asymptotic behavior of the expectations \ $(\EE(Y_t), \EE(X_t))$
 \ as \ $t \to \infty$, \ we introduce a classification of the Heston model
 given by the SDE \eqref{Heston_SDE}.

\begin{Def}\label{Def_criticality}
Let \ $(Y_t, X_t)_{t\in\RR_+}$ \ be the unique strong solution of the SDE
 \eqref{Heston_SDE} satisfying \ $\PP(Y_0 \in \RR_+) = 1$.
\ We call \ $(Y_t, X_t)_{t\in\RR_+}$ \ subcritical, critical or supercritical if
 \ $b \in \RR_{++}$, \ $b = 0$ \ or \ $b \in \RR_{--}$, \ respectively.
\end{Def}

In the sequel \ $\stoch$, \ $\distr$ \ and \ $\as$ \ will denote convergence
 in probability, in distribution and almost surely, respectively.

The following result states the existence of a unique stationary distribution
 and the ergodicity for the process \ $(Y_t)_{t\in\RR_+}$ \ given by the first
 equation in \eqref{Heston_SDE} in the subcritical case, see, e.g., Cox et
 al.\ \cite[Equation (20)]{CoxIngRos}, Li and Ma \cite[Theorem 2.6]{LiMa},
 Theorem 3.1 with \ $\alpha = 2$ \ and Theorem 4.1 in Barczy et al.\
 \cite{BarDorLiPap2}, or Jin et al.\ \cite[Corollaries 5.9 and 6.4]{JinRudTra}.
Only \eqref{ergodic_discrete1} of the following Theorem \ref{Ergodicity} can
 be considered as a slight improvement of the existing results.

\begin{Thm}\label{Ergodicity}
Let \ $a, b, \sigma_1 \in \RR_{++}$.
\ Let \ $(Y_t)_{t\in\RR_+}$ \ be the unique strong solution of the first equation
 of the SDE \eqref{Heston_SDE} satisfying \ $\PP(Y_0 \in \RR_+) = 1$.
\ Then
 \renewcommand{\labelenumi}{{\rm(\roman{enumi})}}
 \begin{enumerate}
  \item
   $Y_t \distr Y_\infty$ \ as \ $t \to \infty$, \ and the distribution of
   \ $Y_\infty$ \ is given by
   \begin{align}\label{Laplace}
    \EE(\ee^{-\lambda Y_\infty})
    = \left(1 + \frac{\sigma_1^2}{2b} \lambda \right)^{-2a/\sigma_1^2} ,
    \qquad \lambda \in \RR_+ ,
   \end{align}
   i.e., \ $Y_\infty$ \ has Gamma distribution with parameters
   \ $2a / \sigma_1^2$ \ and \ $2b / \sigma_1^2$, \ hence
   \begin{align}\label{stac_moment}
     \EE(Y_\infty) = \frac{a}{b}, \qquad
     \EE(Y_\infty^2) = \frac{(2a+\sigma_1^2)a}{2b^2} , \qquad
     \EE(Y_\infty^3) = \frac{(2a+\sigma_1^2)(a+\sigma_1^2)a}{2b^3} .
   \end{align}
 \item
  supposing that the random initial value \ $Y_0$ \ has the same distribution
  as \ $Y_\infty$, \ the process \ $(Y_t)_{t\in\RR_+}$ \ is strictly stationary.
 \item
  for all Borel measurable functions \ $f : \RR \to \RR$ \ such that
  \ $\EE(|f(Y_\infty)|) < \infty$, \ we have
  \begin{equation}\label{ergodic}
   \frac{1}{T} \int_0^T f(Y_s) \, \dd s \as \EE(f(Y_\infty)) \qquad
   \text{as \ $T \to \infty$,}
  \end{equation}
  \begin{equation}\label{ergodic_discrete1}
   \frac{1}{n}\sum_{i=0}^{n-1} f(Y_i) \as \EE(f(Y_\infty)) \qquad \text{as \ $n\to\infty$.}
  \end{equation}
\end{enumerate}
\end{Thm}

\noindent{\bf Proof.}
Based on the references given before the theorem, we only need to show \eqref{ergodic_discrete1}.
By Corollary 2.7 in Jin et al.\ \cite{JinManRudTra}, the tail $\sigma$-field
 \ $\bigcap_{t\in\RR_+} \sigma(Y_s,s\geq t)$ \ of \ $(Y_t)_{t \in \RR_+}$ \ is trivial for any initial distribution,
 i.e., the tail \ $\sigma$-field in question consists of events having probability \ $0$ \ or \ $1$ \
 for any initial distribution on \ $\RR_+$.
\ But since the tail $\sigma$-field of $(Y_t)_{t \in \RR_+}$ is richer than that of $(Y_i)_{i \in \ZZ_+}$,
 the tail \ $\sigma$-field of \ $(Y_i)_{i\in\ZZ_+}$ \ is also trivial for any initial distribution.

Denoting the distribution of \ $Y_0$ \ {and \ $Y_\infty$} \ by \ $\nu$ \ and \ $\mu$, \ respectively, let us introduce the distribution
 \ $\eta := (\mu + \nu)/2$.
\ Let us introduce the following processes: \ $(Z_t)_{t \in \RR_+}$, \ which is
 the pathwise unique strong solution of the first equation in \eqref{Heston_SDE} with initial condition \ $Z_0 = \zeta_0$, \ where \ $\zeta_0$ \ has the distribution
 \ $\mu$; \ and \ $(U_t)_{t \in \RR_+}$, \ which is the pathwise unique strong solution of the same SDE with initial condition \ $U_0 = \xi_0$, \ where \ $\xi_0$ \ has the distribution $\eta$.

We use Birkhoff's ergodic theorem (see, e.g., Theorem 8.4.1 in Dudley
 \cite{Dud}) in the usual setting: the probability space is
 \ $(\RR^{\ZZ_+}, \cB(\RR^{\ZZ_+}), \cL((Z_i)_{i\in\ZZ_+}))$, \ where
 \ $\cL((Z_i)_{i\in\ZZ_+})$ \ denotes the distribution of \ $(Z_i)_{i\in\ZZ_+}$,
 \ and the measure-preserving transformation \ $T$ \ is the shift operator, i.e.,
 \ $T((x_i)_{i\in\ZZ_+}):=(x_{i+1})_{i\in\ZZ_+}$ \ for \ $(x_i)_{i\in\ZZ_+}\in\RR^{\ZZ_+}$ \
 (the measure preservability follows from (ii)).
\ All invariant sets of \ $T$ \ are included in the tail \ $\sigma$-field
 \ of the coordinate mappings \ $\pi_i$, $i\in\ZZ_+$, \ on \ $\RR^{\ZZ_+}$,
 \ since for any invariant set \ $A$ \ we have
 \ $A \in \sigma(\pi_0, \pi_1, \ldots)$, \ but as \ $T^k(A) = A$ \ for all \ $k \in \NN$,
 \ it is also true that \ $A \in \sigma(\pi_k, \pi_{k+1}, \ldots)$ \ for
 all \ $k \in \NN$.
\ This implies that \ $T$ \ is ergodic, since the tail \ $\sigma$-field is trivial.
Hence we can apply the ergodic theorem for the function
\[
g:\RR^{\ZZ_+} \to \RR, \qquad g((x_i)_{i\in\ZZ_+}) := f(x_0) , \qquad (x_i)_{i\in\ZZ_+}\in\RR^{\ZZ_+},
\]
 where \ $f$ \ is given in (iii), to obtain
 \[
   \frac{1}{n}\sum_{i=0}^{n-1} f(x_i) \to \int_{\RR_+} f(x_0)\,\mu(\dd x_0)
   \qquad \text{as \ $n\to\infty$}
 \]
 for almost every \ $(x_i)_{i\in\ZZ_+}\in\RR^{\ZZ_+}$ \ with respect to the measure \ $\cL((Z_i)_{i\in\ZZ_+})$,
 \ and consequently
\begin{equation}\label{Z_convergence}
\frac{1}{n}\sum_{i=0}^{n-1} f(Z_i) \as \EE(f(Y_\infty)) \qquad
 \text{as \ $n \to \infty$,}
\end{equation}
because, clearly, the distribution of \ $Y_\infty$ \ does not depend on the initial distribution. We introduce the following event, which is clearly a tail event of \ $(Z_i)_{i \in \ZZ_+}$ \ and has probability 1 by \eqref{Z_convergence}:
\[
C_Z:=\left\{\omega \in \Omega:
 \text{$\frac{1}{n}\sum_{i=0}^{n-1} f(Z_i(\omega)) \to \EE(f(Y_\infty))$ \ as \ $n \to \infty$}
\right\}.
\]
The events \ $C_Y$ \ and \ $C_U$ \ are defined in a similar way and are clearly tail events of \ $(Y_i)_{i \in \ZZ_+}$ \ and \ $(U_i)_{i \in \ZZ_+}$, \ respectively.
Clearly,
\begin{align*}
\PP(C_U) = \int_{0}^{\infty}\PP(C_U\mid U_0=x) \, \dd \eta(x) &= \frac{1}{2}\int_{0}^{\infty}\PP(C_U \mid U_0=x) \, \dd \mu(x)
 + \frac{1}{2}\int_{0}^{\infty}\PP(C_U \mid U_0=x) \, \dd \nu(x) \\
& \geq \frac{1}{2}\int_{0}^{\infty}\PP(C_U \mid U_0=x) \, \dd \mu(x) = \frac{1}{2}\int_{0}^{\infty}\PP(C_Z \mid Z_0=x) \, \dd \mu(x) \\
& = \frac{1}{2} \PP(C_Z) = \frac{1}{2}.
\end{align*}
Here we used that \ $\PP(C_U \mid U_0=x)= \PP(C_Z \mid Z_0=x)$ \ $\mu$-a.e. \ $x\in\RR_+$, \
 since the conditional probabilities on both sides depend only on the transition
 probability kernel of the CIR process given by the first SDE of \eqref{Heston_SDE} irrespective of the initial distribution.
Further, we note that \ $\PP(C_U \mid U_0=x)$ \ is defined uniquely only \ $\eta$-a.e. \ $x\in\RR_+$, \
 but, by the definition of \ $\eta$, \ this means both \ $\mu$-a.e. \ $x\in\RR_+$, \ and \ $\nu$-a.e. \ $x\in\RR_+$, \
 and similarly \ $\PP(C_Z \mid Z_0=x)$ \ is defined \ $\mu$-a.e. \ $x\in\RR_+$, \ so our equalities are valid.
Thus, we have \ $\PP(C_U) \geq \frac{1}{2}$. \ But since \ $C_U$ \ is a tail event of \ $(U_i)_{i \in \ZZ_+}$, \ its probability must be either 0 or 1 (since the tail \ $\sigma$-field is trivial), hence \ $\PP(C_U)=1$.
\ Hence
 \begin{align*}
  2 = \int_0^\infty \PP(C_U \mid U_0=x) \,\dd \mu(x)
      + \int_0^\infty \PP(C_U \mid U_0=x) \,\dd \nu(x)
    \leq \mu([0,\infty)) +  \nu([0,\infty))
    =2,
 \end{align*}
 yielding that
 \[
     \int_{0}^{\infty}\PP(C_U \mid U_0=x) \, \dd \mu(x) = \int_{0}^{\infty}\PP(C_U \mid U_0=x) \, \dd \nu(x) = 1,
 \]
 and the second equality is exactly \eqref{ergodic_discrete1} after we note that, by the same argument as above,
 \[
 \int_{0}^{\infty}\PP(C_U \mid U_0=x) \, \dd \nu(x)
  = \int_{0}^{\infty}\PP(C_Y \mid Y_0=x) \, \dd \nu(x) = \PP(C_Y).
 \]
With this our proof is complete.
\proofend

In what follows we recall some limit theorems for (local) martingales.
We will use these limit theorems later on for studying the asymptotic
 behaviour of (conditional) least squares estimators for \ $(a, b, \alpha, \beta)$.

First, we recall a strong law of large numbers for discrete time
 square-integrable martingales.

\begin{Thm}{\bf (Shiryaev \cite[Chapter VII, Section 5, Theorem 4]{Shi})}
\label{DDS_stoch_int_discr}
Let \ $\bigl( \Omega, \cF, (\cF_n)_{n\in\NN}, \PP \bigr)$ \ be a filtered
 probability space.
Let \ $(M_n)_{n\in\NN}$ \ be a square-integrable martingale with respect
 to the filtration \ $(\cF_n)_{n\in\NN}$ \ such that \ $\PP(M_0=0)=1$ \ and
 \ $\PP(\lim_{n\to\infty} \langle M\rangle_n = \infty)=1$, \ where
 \ $(\langle M\rangle_n)_{n\in\NN}$ \ denotes the predictable quadratic variation process of
 \ $M$.
\ Then
 \[
   \frac{M_n}{\langle M \rangle_n} \as 0 \qquad \text{as \ $n \to \infty$.}
 \]
\end{Thm}

Next, we recall a martingale central limit theorem in discrete time.

\begin{Thm}{\bf (Jacod and Shiryaev \cite[Chapter VIII, Theorem 3.33]{JSh})}%
\label{Thm_CLT_discrete}
Let \ $\{ (\bM_{n,k},\cF_{n,k}) : k=0,1,\ldots,k_n \}_{n\in\NN}$ \ be a
 sequence of \ $d$-dimensional square-integrable martingales with
 \ $\bM_{n,0} = \bzero$ \ such that there exists some symmetric, positive semi-definite non-random matrix
 \ $\bD \in \RR^{d\times d}$ \ such that
 \[
   \sum_{k=1}^{k_n}
    \EE((\bM_{n,k} - \bM_{n,k-1})
        (\bM_{n,k} - \bM_{n,k-1})^\top \mid \cF_{n,k-1})
   \stoch \bD
   \qquad \text{as \ $n \to \infty$,}
 \]
 and for all \ $\vare \in \RR_{++}$,
 \begin{align}\label{Lindeberg}
  \sum_{k=1}^{k_n}
   \EE( \|\bM_{n,k} - \bM_{n,k-1}\|^2
        \bbone_{\{\|\bM_{n,k}-\bM_{n,k-1}\|\geq\vare\}}
        \mid \cF_{n,k-1})
   \stoch 0 \qquad \text{as \ $n \to \infty$.}
 \end{align}
Then
 \[
   \sum_{k=1}^{k_n}(\bM_{n,k} - \bM_{n,k-1}) = \bM_{n,k_n}
   \distr \cN_d(\bzero, \bD) \qquad \text{as \ $n \to \infty$,}
 \]
 where \ $\cN_d(\bzero, \bD)$ \ denotes a \ $d$-dimensional normal distribution with mean
 vector \ $\bzero$ \ and covariance matrix \ $\bD$.
\end{Thm}

In all the remaining sections, we will consider the subcritical Heston model
 \eqref{Heston_SDE} with a non-random initial value \ $(y_0, x_0) \in \RR_+ \times \RR$.
\ Note that the augmented filtration \ $(\cF_t)_{t\in\RR_+}$ \ corresponding to \ $(W_t,B_t)_{t\in\RR_+}$
 \ and the initial value \ $(y_0,x_0)\in\RR_+\times\RR$, \ in fact, does not depend on \ $(y_0,x_0)$.

\section{CLSE based on discrete time observations}

Using \eqref{cond_exp_discrete_Y} and \eqref{cond_exp_discrete_X}, by an easy
 calculation, for all \ $i \in \NN$,
 \begin{equation}\label{cond_exp_discrete}
  \EE\left( \begin{bmatrix}
             Y_i \\
             X_i
            \end{bmatrix} \bigg| \, \cF_{i-1} \right)
  = \begin{bmatrix}
     \ee^{-b} & 0 \\
     - \beta \int_0^1 \ee^{-bu} \, \dd u & 1 \\
    \end{bmatrix}
    \begin{bmatrix}
     Y_{i-1} \\
     X_{i-1} \\
    \end{bmatrix}
    + \begin{bmatrix}
       \int_0^1 \ee^{-bu} \, \dd u & 0 \\
       - \beta \int_0^1 \left(\int_0^u \ee^{-bv} \, \dd v \right) \dd u & 1 \\
      \end{bmatrix}
      \begin{bmatrix}
       a \\
       \alpha \\
      \end{bmatrix}.
 \end{equation}
Using that \ $\sigma(X_1,Y_1,\ldots,X_{i-1},Y_{i-1}) \subseteq \cF_{i-1}$, $i\in\NN$, \
 by tower rule for conditional expectations, we have
 \begin{align*}
    &\EE\left( \begin{bmatrix}
               Y_i \\
               X_i
              \end{bmatrix} \Bigg| \, \sigma(X_1,Y_1,\ldots,X_{i-1},Y_{i-1}) \right)
    = \EE\left(
       \EE\left( \begin{bmatrix}
                 Y_i \\
                 X_i
                \end{bmatrix} \Bigg| \, \cF_{i-1} \right)
                \Bigg| \, \sigma(X_1,Y_1,\ldots,X_{i-1},Y_{i-1})
                \right) \\
    &\qquad  = \begin{bmatrix}
     \ee^{-b} & 0 \\
     - \beta \int_0^1 \ee^{-bu} \, \dd u & 1 \\
    \end{bmatrix}
    \begin{bmatrix}
     Y_{i-1} \\
     X_{i-1} \\
    \end{bmatrix}
    + \begin{bmatrix}
       \int_0^1 \ee^{-bu} \, \dd u & 0 \\
       - \beta \int_0^1 \left(\int_0^u \ee^{-bv} \, \dd v \right) \dd u & 1 \\
      \end{bmatrix}
      \begin{bmatrix}
       a \\
       \alpha \\
      \end{bmatrix}, \qquad i\in\NN,
 \end{align*}
 and hence a CLSE of \ $(a, b, \alpha, \beta)$ \ based on discrete time observations
 \ $(Y_i, X_i)_{i\in\{1,\ldots,n\}}$ \ could be obtained by solving the extremum
 problem
 \begin{align}\label{CLSE_def}
  \argmin_{(a,b,\alpha,\beta)\in\RR^4}
      \sum_{i=1}^n
       \left[ (Y_i - d Y_{i-1} - c)^2
              + (X_i - X_{i-1} - \gamma - \delta Y_{i-1})^2 \right] ,
 \end{align}
 where
 \begin{equation}\label{c_d_gamma_delta}
 \begin{split}
  d &:= d(b) := \ee^{-b} , \qquad
  c := c(a,b) := a \int_0^1 \ee^{-bu} \, \dd u , \\
  \delta &:= \delta(b, \beta) := -\beta \int_0^1 \ee^{-bu} \, \dd u , \qquad
  \gamma := \gamma(a, b, \alpha, \beta)
         := \alpha
            - a \beta \int_0^1 \left(\int_0^u \ee^{-bv} \, \dd v \right) \dd u.
 \end{split}
 \end{equation}
 First, we determine the CLSE of \ $(c, d, \gamma, \delta)$ \ by minimizing the
 sum on the right hand side of \eqref{CLSE_def} with respect to
 \ $(c, d, \gamma, \delta) \in \RR^4$.

We get
\begin{align}\label{CLSEcdgammadelta_discrete}
  \begin{bmatrix}
   \hc_n^{\CLSED} \\
   \hd_n^{\CLSED} \\
   \hgamma_n^{\CLSED} \\
   \hdelta_n^{\CLSED}
  \end{bmatrix}
  &= \left( \bI_2 \otimes \begin{bmatrix}
                           n &   \sum_{i=1}^n Y_{i-1} \\
                           \sum_{i=1}^n Y_{i-1} & \sum_{i=1}^n Y_{i-1}^2
                          \end{bmatrix}^{-1} \right)
     \begin{bmatrix}
      \sum_{i=1}^n Y_{i} \\
      \sum_{i=1}^n Y_i Y_{i-1} \\
      X_n - x_0 \\
      \sum_{i=1}^n (X_i - X_{i-1}) Y_{i-1}
  \end{bmatrix}
 \end{align}
 provided that \ $n \sum_{i=1}^n Y_{i-1}^2 > \left(\sum_{i=1}^n Y_{i-1}\right)^2$,
 \ where \ $\otimes$ \ denotes Kronecker product of matrices.
Indeed, with the notation
 \[
   f(c, d, \gamma, \delta)
   := \sum_{i=1}^n
       \left[ (Y_i - d Y_{i-1} - c)^2
              + (X_i - X_{i-1} - \gamma - \delta Y_{i-1})^2 \right] ,
   \qquad (c, d, \gamma, \delta) \in \RR^4 ,
 \]
 we have
 \begin{align*}
  &\frac{\partial f}{\partial c}(c, d, \gamma, \delta)
   = -2\sum_{i=1}^n (Y_i - dY_{i-1} - c),\\
  &\frac{\partial f}{\partial d}(c, d, \gamma, \delta)
   = -2\sum_{i=1}^n Y_{i-1}(Y_i - dY_{i-1} - c),\\
  &\frac{\partial f}{\partial \gamma}(c, d, \gamma, \delta)
   = -2\sum_{i=1}^n (X_i - X_{i-1} - \gamma - \delta Y_{i-1}),\\
  &\frac{\partial f}{\partial \delta}(c, d, \gamma, \delta)
   = -2\sum_{i=1}^n Y_{i-1}(X_i - X_{i-1} - \gamma - \delta Y_{i-1}) .
 \end{align*}
Hence the system of equations consisting of the first order partial derivates
 of \ $f$ \ being equal to \ $0$ \ takes the form
 \[
   \left( \bI_2 \otimes \begin{bmatrix}
                         n & \sum_{i=1}^n Y_{i-1} \\
                         \sum_{i=1}^n Y_{i-1} & \sum_{i=1}^n Y_{i-1}^2
                        \end{bmatrix} \right)
   \begin{bmatrix}
    c \\
    d \\
    \gamma \\
    \delta
   \end{bmatrix}
   = \begin{bmatrix}
      \sum_{i=1}^n Y_i \\
      \sum_{i=1}^n Y_{i-1} Y_i \\
      X_n - x_0 \\
      \sum_{i=1}^n (X_i - X_{i-1}) Y_{i-1}
     \end{bmatrix} .
 \]
This implies \eqref{CLSEcdgammadelta_discrete}, since the $4\times 4$-matrix
 consisting of the second order partial derivatives of \ $f$ \ having the form
 \[
   2 \bI_2 \otimes \begin{bmatrix}
                  n & \sum_{i=1}^n Y_{i-1} \\
                  \sum_{i=1}^n Y_{i-1} & \sum_{i=1}^n Y_{i-1}^2
                 \end{bmatrix}
 \]
 is positive definite provided that
 \ $n \sum_{i=1}^n Y_{i-1}^2 > \left(\sum_{i=1}^n Y_{i-1}\right)^2$.
\ In fact, it turned out that for the calculation of the CLSE of
 \ $(c, d, \gamma, \delta)$, \ one does not need to know the values of the
 parameters \ $\sigma_1, \sigma_2$ \ and \ $\varrho$.

The next lemma assures the unique existence of the CLSE of
 \ $(c,d,\gamma,\delta)$ \ based on
 discrete time observations.

\begin{Lem}\label{LEMMA_LSE_exist_discrete}
If \ $a \in \RR_{++}$, \ $b \in \RR$, \ $\sigma_1\in\RR_{++}$, \ and
 \ $Y_0 = y_0  \in \RR_+$, \
 then for all \ $n\geq 2$, \ $n\in\NN$, \ we have
 \[
  \PP\left( n\sum_{i=1}^n Y_{i-1}^2 > \left(\sum_{i=1}^n Y_{i-1}\right)^2\right)=1,
 \]
 and hence, supposing also that \ $\alpha,\beta\in\RR$, $\sigma_2\in\RR_{++}$, \ $\varrho \in(-1,1)$,
 there exists a unique CLSE \ $(\hc_n^{\CLSED}, \hd_n^{\CLSED}, \hgamma_n^{\CLSED}, \hdelta_n^{\CLSED})$
 \ of \ $(c,d,\gamma,\delta)$ \ which has the form given in \eqref{CLSEcdgammadelta_discrete}.
\end{Lem}

\noindent{\bf Proof.}
By an easy calculation,
 \[
   n\sum_{i=1}^n Y_{i-1}^2 - \left(\sum_{i=1}^n Y_{i-1}\right)^2
     = n \sum_{i=1}^n \left(  Y_{i-1} - \frac{1}{n}\sum_{j=1}^n Y_{j-1} \right)^2
     \geq 0,
 \]
 and equality holds if and only if
 \[
    Y_{i-1} = \frac{1}{n}\sum_{j=1}^n Y_{j-1},
    \qquad i=1,\ldots,n
    \qquad  \Longleftrightarrow \qquad
    Y_0=Y_1=\cdots=Y_{n-1}.
 \]
Then, for all \ $n\geq 2$,
 \[
   \PP(Y_0 = Y_1 = \cdots = Y_{n-1})
       \leq \PP(Y_0=Y_1) = \PP(Y_1=y_0)=0,
 \]
 since the law of \ $Y_1$ \ is absolutely continuous, see, e.g., Cox et al. \cite[formula 18]{CoxIngRos}.
\proofend

Note that Lemma \ref{LEMMA_LSE_exist_discrete} is valid for all \ $b\in\RR$, \ i.e.,
 not only for the subcritical Heston model.

Next, we describe the asymptotic behaviour of the CLSE of \ $(c,d,\gamma,\delta)$.

\begin{Thm}\label{c_d_normal}
If \ $a, b \in \RR_{++}$, \ $\alpha, \beta \in \RR$,
 \ $\sigma_1, \sigma_2 \in \RR_{++}$, \ $\varrho \in (-1, 1)$ \ and
 \ $(Y_0, X_0) = (y_0, x_0) \in \RR_{++} \times \RR$, \ then the CLSE
 \ $(\hc_n^{\CLSED}, \hd_n^{\CLSED}, \hgamma_n^{\CLSED}, \hdelta_n^{\CLSED})$
 \ of \ $(c, d, \gamma, \delta)$ \ given in \eqref{CLSEcdgammadelta_discrete}
 is strongly consistent and asymptotically normal, i.e.,
 \[
   (\hc_n^{\CLSED}, \hd_n^{\CLSED}, \hgamma_n^{\CLSED}, \hdelta_n^{\CLSED})
   \as (c, d, \gamma, \delta) \qquad \text{as \ $n \to \infty$,}
 \]
 and
 \[
   \sqrt{n}
   \begin{bmatrix}
    \hc_n^{\CLSED} - c \\
    \hd_n^{\CLSED} - d \\
    \hgamma_n^{\CLSED} - \gamma \\
    \hdelta_n^{\CLSED} - \delta
   \end{bmatrix}
   \distr
   \cN_4\left(\bzero, \bE  \right)
   \qquad \text{as \ $n \to \infty$,}
 \]
 with some explicitly given symmetric, positive definite matrix
 \ $\bE \in \RR^{2\times 2}$ \ given in \eqref{matrixE}.
\end{Thm}

{\bf \noindent Proof.}
By \eqref{CLSEcdgammadelta_discrete}, we get
 \begin{equation}\label{help3}
  \begin{split}
   \begin{bmatrix}
    \hc_n^{\CLSED} \\
    \hd_n^{\CLSED}
   \end{bmatrix}
   &=\left( \sum_{i=1}^n
 	     \begin{bmatrix} 1 \\ Y_{i-1} \end{bmatrix}
 	     \begin{bmatrix} 1 \\ Y_{i-1} \end{bmatrix}^\top \right)^{-1}
     \left( \sum_{i=1}^{n}
             \begin{bmatrix} 1 \\ Y_{i-1} \end{bmatrix} Y_i \right) \\
   &=\left( \sum_{i=1}^n
             \begin{bmatrix} 1 \\ Y_{i-1} \end{bmatrix}
   	     \begin{bmatrix} 1 \\ Y_{i-1} \end{bmatrix}^\top \right)^{-1}
     \left( \sum_{i=1}^n
       	     \begin{bmatrix} 1 \\ Y_{i-1} \end{bmatrix}
       	     \begin{bmatrix} 1 \\ Y_{i-1} \end{bmatrix}^\top \right)
     \begin{bmatrix} c \\ d \end{bmatrix} \\
   &\quad
     + \left( \sum_{i=1}^n
               \begin{bmatrix} 1 \\ Y_{i-1} \end{bmatrix}
               \begin{bmatrix} 1 \\ Y_{i-1} \end{bmatrix}^\top \right)^{-1}
       \sum_{i=1}^{n}
        \begin{bmatrix} 1 \\ Y_{i-1} \end{bmatrix} (Y_i - c - d Y_{i-1}) \\
   &= \begin{bmatrix} c \\ d \end{bmatrix}
      + \left( \frac{1}{n}
               \sum_{i=1}^n
              	\begin{bmatrix} 1 \\ Y_{i-1} \end{bmatrix}
              	\begin{bmatrix} 1 \\ Y_{i-1} \end{bmatrix}^\top \right)^{-1}
        \frac{1}{n}
        \sum_{i=1}^{n} \begin{bmatrix} 1 \\ Y_{i-1} \end{bmatrix} \vare_i ,
  \end{split}
 \end{equation}
 where \ $\vare_i := Y_i - c - d Y_{i-1}$, \ $i \in \NN$, \ provided that
 \ $n \sum_{i=1}^n Y_{i-1}^2 > \left(\sum_{i=1}^n Y_{i-1}\right)^2$.
\ By \eqref{cond_exp_discrete} and \eqref{c_d_gamma_delta}, \ $\EE(Y_i \mid \cF_{i-1}) = d Y_{i-1} + c$,
 \ $i \in \NN$, \ and hence \ $(\vare_i)_{i\in\NN}$ \ is a sequence of
 martingale differences with respect to the filtration \ $(\cF_i)_{i\in\ZZ_+}$.
\ By \eqref{Solutions}, we have
 \begin{align*}
  Y_i &= \ee^{-b} Y_{i-1} + a \int_{i-1}^i \ee^{-b(i-u)} \, \dd u
         + \sigma_1 \int_{i-1}^i \ee^{-b(i-u)} \sqrt{Y_u} \, \dd W_u \\
      &= d Y_{i-1} + c
         + \sigma_1 \int_{i-1}^i \ee^{-b(i-u)} \sqrt{Y_u} \, \dd W_u ,
  \qquad i \in \NN ,
 \end{align*}
 hence, by Proposition 3.2.10 in Karatzas and Shreve \cite{KarShr}
 and \eqref{cond_exp_discrete_Y}, we have
 \begin{align*}
  \EE(\vare_i^2\mid\cF_{i-1})
  &= \sigma_1^2
     \EE\biggl(\left( \int_{i-1}^i \ee^{-b(i-u)} \sqrt{Y_u} \, \dd W_u \right)^2
               \, \bigg| \, \cF_{i-1} \biggr)
   = \sigma_1^2 \int_{i-1}^i \ee^{-2b(i-u)} \EE(Y_u \mid \cF_{i-1}) \, \dd u \\
  &= \sigma_1^2 \int_{i-1}^i \ee^{-2b(i-u)} \ee^{-b(u-i+1)} Y_{i-1} \, \dd u
     + \sigma_1^2 \int_{i-1}^{i}\ee^{-2b(i-u)}
       a \int_{i-1}^{u} \ee^{-b(u-v)} \, \dd v \, \dd u \\
  &= \sigma_1^2 Y_{i-1} \int_0^1 \ee^{-b(2-v)} \, \dd v
     + \sigma_1^2 a \int_0^1 \int_0^u \ee^{-b(2-v-u)} \, \dd v \, \dd u
   =: C_1 Y_{i-1} + C_2 .
 \end{align*}
Now we apply Theorem \ref{DDS_stoch_int_discr} to the square-integrable martingale
 \ $M_n^{(c)} := \sum_{i=1}^{n} \vare_i$, \ $n \in \NN$, \ which has predictable quadratic
 variation process
 \ $\langle M^{(c)} \rangle_n = \sum_{i=1}^n \EE(\vare_i^2 \mid \cF_{i-1})
    = C_1 \sum_{i=1}^n Y_{i-1} + C_2 n$,
 \ $n \in \NN$,
 \ see, e.g., Shiryaev \cite[Chapter VII, Section 1, formula (15)]{Shi}.
\ By \eqref{stac_moment} and \eqref{ergodic_discrete1},
 \[
   \frac{\langle M^{(c)} \rangle_n}{n} \as C_1 \EE(Y_\infty) + C_2 \qquad
   \text{as \ $n \to \infty$,}
 \]
 and since \ $C_1,C_2\in\RR_{++}$, \ $\langle M^{(c)}\rangle_n \as \infty$ \ as
 \ $n \to \infty$.
\ Hence, by Theorem \ref{DDS_stoch_int_discr},
 \begin{equation}\label{vare_disappear}
   \frac{1}{n} \sum_{i=1}^n \vare_i
   = \frac{M_n^{(c)}}{\langle M^{(c)} \rangle_n}
     \frac{\langle M^{(c)} \rangle_n}{n}
   \as 0 \cdot (C_1 \EE(Y_{\infty}) + C_2) = 0
   \qquad \text{as \ $n\to\infty$.}
 \end{equation}
Similarly,
 \begin{align*}
  \EE(Y_{i-1}^2 \vare_i^2 \mid \cF_{i-1})
  = Y_{i-1}^2 \EE(\vare_i^2\mid \cF_{i-1})
  = C_1 Y_{i-1}^3 + C_2Y_{i-1}^2,
  \qquad i \in \NN ,
 \end{align*}
 and, by essentially the same reasoning as before,
 \ $\frac{1}{n} \sum_{i=1}^n Y_{i-1} \vare_i \as 0$ \ as \ $n \to \infty$.
\ By \eqref{stac_moment} and \eqref{ergodic_discrete1},
 \begin{align}\label{help4}
  \left( \frac{1}{n}
         \sum_{i=1}^n
          \begin{bmatrix} 1 \\ Y_{i-1} \end{bmatrix}
          \begin{bmatrix} 1 \\ Y_{i-1} \end{bmatrix}^\top \right)^{-1}
  = \begin{bmatrix}
     1 & \frac{1}{n} \sum_{i=1}^n Y_{i-1} \\
     \frac{1}{n} \sum_{i=1}^n Y_{i-1} & \frac{1}{n} \sum_{i=1}^n Y_{i-1}^2
    \end{bmatrix}^{-1}
  \as \begin{bmatrix}
       1 & \EE(Y_\infty) \\
       \EE(Y_\infty) & \EE(Y_\infty^2)
      \end{bmatrix}^{-1}
 \end{align}
 as \ $n \to \infty$, where we used that
 \ $\EE(Y_\infty^2) - (\EE(Y_\infty))^2 = \frac{a\sigma_1^2}{2b^2} \in \RR_{++}$,
 \ and consequently, the limit is indeed non-singular.
Thus, by \eqref{help3}, \ $(\hc_n^{\CLSED}, \hd_n^{\CLSED}) \as (c, d)$ \ as
 \ $n \to \infty$.

Further, by \eqref{CLSEcdgammadelta_discrete},
 \begin{equation}\label{help6}
  \begin{split}
   \begin{bmatrix} \hgamma_n^{\CLSED} \\ \hdelta_n^{\CLSED} \end{bmatrix}
   &= \left( \sum_{i=1}^n
              \begin{bmatrix} 1 \\ Y_{i-1} \end{bmatrix}
 	      \begin{bmatrix} 1 \\ Y_{i-1} \end{bmatrix}^\top \right)^{-1}
      \left( \sum_{i=1}^n
   	      \begin{bmatrix} 1 \\ Y_{i-1} \end{bmatrix}
              (X_i - X_{i-1}) \right) \\
   &= \left( \sum_{i=1}^n
   	      \begin{bmatrix} 1 \\ Y_{i-1} \end{bmatrix}
   	      \begin{bmatrix} 1 \\ Y_{i-1} \end{bmatrix}^\top \right)^{-1}
      \left( \sum_{i=1}^n
       	      \begin{bmatrix} 1 \\ Y_{i-1} \end{bmatrix}
       	      \begin{bmatrix} 1 \\ Y_{i-1} \end{bmatrix}^\top \right)
      \begin{bmatrix}
      \gamma \\ \delta
      \end{bmatrix} \\
   &\quad
      + \left( \sum_{i=1}^n
         	\begin{bmatrix} 1 \\ Y_{i-1} \end{bmatrix}
         	\begin{bmatrix} 1 \\ Y_{i-1} \end{bmatrix}^\top \right)^{-1}
        \sum_{i=1}^n
         \begin{bmatrix} 1 \\ Y_{i-1} \end{bmatrix}
         (X_i - X_{i-1} - \gamma - \delta Y_{i-1}) \\
   &= \begin{bmatrix} \gamma \\ \delta \end{bmatrix}
      + \left( \frac{1}{n}
               \sum_{i=1}^n
              	\begin{bmatrix} 1 \\ Y_{i-1} \end{bmatrix}
              	\begin{bmatrix} 1 \\ Y_{i-1} \end{bmatrix}^\top \right)^{-1}
        \frac{1}{n}
        \sum_{i=1}^n
         \begin{bmatrix} 1 \\ Y_{i-1} \end{bmatrix}
         \eta_i ,
  \end{split}
 \end{equation}
 where \ $\eta_i := X_i - X_{i-1} - \gamma - \delta Y_{i-1}$, \ $i \in \NN$,
 \ provided that
 \ $n \sum_{i=1}^n Y_{i-1}^2 > \left(\sum_{i=1}^n Y_{i-1}\right)^2$.
\ By \eqref{cond_exp_discrete} and \eqref{c_d_gamma_delta},
 \ $\EE(X_i \mid \cF_{i-1}) = X_{i-1} + \delta Y_{i-1} + \gamma$, \ $i \in \NN$,
 \ and hence \ $(\eta_i)_{i\in\NN}$ \ is a sequence of martingale differences
 with respect to the filtration $(\cF_i)_{i\in\ZZ_+}$.
\ By \eqref{Solutions} and \eqref{cond_exp_discrete_Y}, with the
 notation \ $\widetilde W_t := \varrho W_t + \sqrt{1 - \varrho^2} B_t$,
 \ $t \in \RR_+$, \ we compute
 \begin{align*}
  X_i - X_{i-1} &= \int_{i-1}^i (\alpha - \beta Y_u) \, \dd u
         + \sigma_2 \int_{i-1}^i \sqrt{Y_u} \, \dd \tW_u
       = \alpha - \beta \int_{i-1}^i Y_u \, \dd u
         + \sigma_2 \int_{i-1}^i \sqrt{Y_u} \, \dd \tW_u \\
      &= \alpha
         - \beta
           \int_{i-1}^i
            \left( \ee^{-b(u-(i-1))} Y_{i-1}
                   + a \int_{i-1}^u \ee^{-b(u-v)} \, \dd v
                   + \sigma_1
                     \int_{i-1}^u \ee^{-b(u-v)} \sqrt{Y_v} \, \dd W_v \right)
            \dd u \\
      &\quad
         + \sigma_2 \int_{i-1}^i \sqrt{Y_u} \, \dd \tW_u \\
      &= \alpha - \beta Y_{i-1} \int_{i-1}^i \ee^{-b(u-i+1)} \, \dd u
         - a \beta
           \int_{i-1}^i \left( \int_{i-1}^u \ee^{-b(u-v)} \, \dd v \right) \dd u \\
      &\quad
         - \beta \sigma_1
           \int_{i-1}^i
            \left( \int_{i-1}^u \ee^{-b(u-v)} \sqrt{Y_v} \, \dd W_v \right) \dd u
         + \sigma_2 \int_{i-1}^i \sqrt{Y_u} \, \dd \tW_u
 \end{align*}
 \begin{align*}
      &= \alpha - \beta Y_{i-1} \int_{0}^1 \ee^{-bv} \, \dd v
         - a \beta
           \int_{0}^1 \left( \int_{0}^u \ee^{-bv} \,\dd v \right) \dd u \\
      &\quad
         - \beta \sigma_1
           \int_{i-1}^i
            \left( \int_{i-1}^u \ee^{-b(u-v)} \sqrt{Y_v} \, \dd W_v \right) \dd u
         + \sigma_2 \int_{i-1}^i \sqrt{Y_u} \, \dd \tW_u \\
      &= \delta Y_{i-1} + \gamma
         - \beta \sigma_1
           \int_{i-1}^i
            \left( \int_{i-1}^u \ee^{-b(u-v)} \sqrt{Y_v} \, \dd W_v \right) \dd u
         + \sigma_2 \int_{i-1}^i \sqrt{Y_u} \, \dd \tW_u ,
 \end{align*}
 and consequently,
 \begin{align*}
  \EE(\eta_i^2 \mid \cF_{i-1})
  &= \beta^2 \sigma_1^2
     \EE\biggl[ \biggl( \int_{i-1}^i \int_{i-1}^u
                         \ee^{-b(u-v)} \sqrt{Y_v} \, \dd W_v\, \dd u \biggr)^2
                \,\bigg|\, \cF_{i-1} \biggr]
     + \sigma_2^2
       \EE\biggl[ \biggl( \int_{i-1}^i \sqrt{Y_u} \, \dd \tW_u \biggr)^2
                  \,\bigg|\, \cF_{i-1} \biggr] \\
  &\quad
     - 2 \beta \sigma_1 \sigma_2
       \EE\biggl[ \bigg( \int_{i-1}^i \int_{i-1}^u
                          \ee^{-b(u-v)} \sqrt{Y_v} \, \dd W_v\, \dd u \biggr)
                  \biggl( \varrho \int_{i-1}^i \sqrt{Y_u} \, \dd W_u \biggr)
                  \,\bigg|\, \cF_{i-1} \biggr] \\
  &\quad
     - 2 \beta \sigma_1 \sigma_2
       \EE\biggl[ \biggl( \int_{i-1}^i \int_{i-1}^u
                           \ee^{-b(u-v)} \sqrt{Y_v} \, \dd W_v \, \dd u \biggr)
                  \biggl( \sqrt{1 - \varrho^2}
                          \int_{i-1}^i \sqrt{Y_u} \, \dd B_u \biggr)
                  \,\bigg|\, \cF_{i-1} \biggr] .
 \end{align*}
We use Equation (3.2.23) from Karatzas and Shreve \cite{KarShr} to the first,
 second
 and third terms, and Proposition 3.2.17 from Karatzas and Shreve
 \cite{KarShr} to the fourth term (together with the independence of \ $W$ \ and \ $B$):
 \begin{align*}
  \EE(\eta_i^2 | \cF_{i-1})
   &= \beta^2 \sigma_1^2 \int_{i-1}^i \int_{i-1}^i
       \EE \left(  \int_{i-1}^u \ee^{-b(u-w)}\sqrt{Y_w}\,\dd W_w
                  \int_{i-1}^v \ee^{-b(v-w)}\sqrt{Y_w}\,\dd W_w
                  \;\Big\vert\; \cF_{i-1}
                   \right)\dd v\,\dd u\\
    &\quad + \sigma_2^2 \int_{i-1}^i \EE(Y_u \mid \cF_{i-1}) \, \dd u \\
    &\quad - 2 \beta \sigma_1 \sigma_2 \varrho
        \int_{i-1}^i \EE\left( \int_{i-1}^u \ee^{-b(u-w)}\sqrt{Y_w}\,\dd W_w
                                 \int_{i-1}^i \sqrt{Y_w}\,\dd W_w \;\Big\vert\; \cF_{i-1} \right)\dd u - 0\\
  &= \beta^2 \sigma_1^2
     \int_{i-1}^i \int_{i-1}^i \int_{i-1}^{u \wedge v}
      \ee^{-b(u+v-2w)} \EE(Y_w \mid \cF_{i-1}) \, \dd w \, \dd u \, \dd v
     + \sigma_2^2 \int_{i-1}^i \EE(Y_u \mid \cF_{i-1}) \, \dd u \\
  &\quad
     - 2 \beta \sigma_1 \sigma_2 \varrho
       \int_{i-1}^i \int_{i-1}^u \ee^{-b(u-v)}
        \EE(Y_v \mid \cF_{i-1}) \, \dd v \, \dd u.
 \end{align*}
Using again \eqref{cond_exp_discrete_Y}, we get
 \begin{align*}
  &\EE(\eta_i^2 | \cF_{i-1})
   = \beta^2 \sigma_1^2 Y_{i-1}
     \int_{i-1}^i \int_{i-1}^i \int_{i-1}^{u \wedge v}
      \ee^{-b(u+v-w-(i-1))} \, \dd w \, \dd v \, \dd u \\
  &\phantom{\EE(\eta_i^2 | \cF_{i-1})\quad}
     + a \beta^2 \sigma_1^2
       \int_{i-1}^i \int_{i-1}^i \int_{i-1}^{u \wedge v} \int_{i-1}^w
        \ee^{-b(u+v-w-z)} \, \dd z \, \dd w \, \dd v \, \dd u
     + \sigma_2^2 Y_{i-1} \int_{i-1}^i e^{-b(u-(i-1))} \, \dd u \\
  &\phantom{\EE(\eta_i^2 | \cF_{i-1})\quad}
     + a \sigma_2^2 \int_{i-1}^i \int_{i-1}^u \ee^{-b(u-v)} \, \dd v \, \dd u
     - 2 \beta \sigma_1 \sigma_2 \varrho Y_{i-1}
       \int_{i-1}^i \int_{i-1}^u \ee^{-b(u-(i-1))} \, \dd v \, \dd u \\
  &\phantom{\EE(\eta_i^2 | \cF_{i-1})\quad}
     - 2 a \beta \sigma_1 \sigma_2 \varrho
       \int_{i-1}^i \int_{i-1}^u \int_{i-1}^v
        \ee^{-b(u-w)} \, \dd w \, \dd v \, \dd u
\end{align*}
 \begin{align*}
  &= \biggl( \beta^2 \sigma_1^2
             \int_{0}^1 \int_{0}^1 \int_{0}^{u' \wedge v'}
              \!\!\ee^{-b(u'+v'-w')} \, \dd w' \, \dd v' \, \dd u'
             - 2 \beta \sigma_1 \sigma_2 \varrho
               \int_{0}^1 \int_{0}^{u'} \ee^{-bu'} \, \dd v' \, \dd u'
             + \sigma_2^2 \int_{0}^1 \!\!\ee^{-bu'} \, \dd u' \biggr) Y_{i-1} \\
  &\quad
     + a \beta^2 \sigma_1^2
       \int_{0}^1 \int_{0}^1 \int_{0}^{u' \wedge v'} \int_{0}^{w'}
       \!\! \ee^{-b(u'+v'-w'-z')} \, \dd z' \, \dd w' \, \dd v' \, \dd u' \\
  &\quad
     + a \sigma_2^2 \int_{0}^1 \int_{0}^{u'} \ee^{-b(u'-v')} \, \dd v' \, \dd u'
     - 2 a \beta \sigma_1 \sigma_2 \varrho
       \int_{0}^1 \int_{0}^{u'} \int_{0}^{v'}
        \ee^{-b(u'-w')} \, \dd w' \, \dd v' \, \dd u'
   =: C_3 Y_{i-1} + C_4 .
\end{align*}
Now we apply Theorem \ref{DDS_stoch_int_discr} to the square-integrable martingale
 \ $M_n^{(\gamma)} := \sum_{i=1}^n \eta_i$, \ $n \in \NN$, \ which has predictable quadratic
 variation process
 \ $\langle M^{(\gamma)} \rangle_n = \sum_{i=1}^n \EE(\eta_i^2 \mid \cF_{i-1})
    = C_3 \sum_{i=1}^n Y_{i-1} + C_4 n$,
 \ $n \in \NN$.
\ By \eqref{ergodic_discrete1},
 \begin{align}\label{help5}
  \frac{\langle M^{(\gamma)} \rangle_n}{n} \as C_3 \EE(Y_\infty) + C_4 \qquad
  \text{as \ $n \to \infty$.}
 \end{align}
Note that \ $C_3 \geq 0$ \ and \ $C_4 \geq 0$, \ since
 \ $\EE(\eta_1^2 \mid \cF_0) = C_3 y_0 + C_4 \geq 0$ \ for all
 \ $y_0 \in \RR_+$.
\ By setting \ $y_0 = 0$, \ we can see that \ $C_4 \geq 0$, \ and
 then, by taking the limit \ $y_0\to\infty$ \ on the right-hand side of the inequality
 \ $C_3\geq -\frac{C_4}{y_0}$, $y_0>0$, \ we get \ $C_3 \geq 0$ \ as well.
Note also that \ $\langle M^{(\gamma)} \rangle_n \as \infty$ \ as
 \ $n \to \infty$ \ provided that \ $C_3 + C_4 > 0$.
\ If \ $C_3 = 0$ \ and \ $C_4 = 0$, \ then \ $\EE(\eta_i^2 \mid \cF_{i-1}) = 0$,
 \ $i \in \NN$, \ and consequently \ $\EE(\eta_i^2)=0$, \ $i\in\NN$, \ and, since \ $\EE(\eta_i) = 0$, \ $i \in \NN$, \ we have
 \ $\PP(\eta_i = 0) = 1$, \ $i \in \NN$, \ implying that
 \ $\PP(\sum_{i=1}^n \eta_i = 0) = 1$ \ and
 \ $\PP(\sum_{i=1}^n Y_{i-1} \eta_i = 0) = 1$, \ $n \in \NN$, \ i.e., in this
 case, by \eqref{help6},
 \ $(\hgamma_n^{\CLSED}, \hdelta_n^{\CLSED}) = (\gamma, \delta)$, \ $n \in \NN$,
 \ almost surely.
If \ $C_3 + C_4 > 0$, \ then, by Theorem  \ref{DDS_stoch_int_discr},
 \begin{equation}\label{eta_disappear}
   \frac{1}{n} \sum_{i=1}^n \eta_i
   = \frac{M_n^{(\gamma)}}{\langle M^{(\gamma)} \rangle_n}
     \frac{\langle M^{(\gamma)} \rangle_n}{n}
   \as 0 \cdot (C_3 \EE(Y_{\infty}) + C_4)
   = 0
   \qquad \text{as \ $n \to \infty$.}
 \end{equation}
Similarly,
 \begin{align*}
  \EE(Y_{i-1}^2 \eta_i^2 \mid \cF_{i-1})
  = Y_{i-1}^2 \EE(\eta_i^2\mid \cF_{i-1})
  = C_3 Y_{i-1}^3 + C_4Y_{i-1}^2 ,
  \qquad i \in \NN ,
 \end{align*}
 and, by essentially the same reasoning as before,
 \ $\frac{1}{n} \sum_{i=1}^n Y_{i-1} \eta_i \as 0$ \ as \ $n \to \infty$ \ (in the case \ $C_3+C_4>0$).
\ Using \eqref{help4} and \eqref{help6}, we have
 \ $(\hgamma_n^{\CLSED}, \hdelta_n^{\CLSED}) \as (\gamma, \delta)$ \ as
 \ $n \to \infty$.

Since the intersection of two events having probability 1 is an event having
 probability 1, we get
 \ $(\hc_n^{\CLSED}, \hd_n^{\CLSED}, \hgamma_n^{\CLSED}, \hdelta_n^{\CLSED})
    \as (c, d, \gamma, \delta)$
 \ as \ $n \to \infty$, \ as desired.

Next, we turn to prove that the CLSE of \ $(c, d, \gamma, \delta)$ \ is
 asymptotically normal.
First, using \eqref{help3} and \eqref{help6}, we can write
 \begin{align}\label{help7}
  \sqrt{n}
  \begin{bmatrix}
   \hc_n^{\CLSED} - c \\
   \hd_n^{\CLSED} - d \\
   \hgamma_n^{\CLSED} - \gamma \\
   \hdelta_n^{\CLSED} - \delta
  \end{bmatrix}
  = \left( \bI_2
           \otimes
           \left( n^{-1}
                  \sum_{i=1}^n
                   \begin{bmatrix} 1 \\ Y_{i-1} \end{bmatrix}
                   \begin{bmatrix} 1 \\ Y_{i-1} \end{bmatrix}^\top
           \right)^{-1} \right)
    n^{-1/2}
    \sum_{i=1}^n
     \begin{bmatrix} \vare_i \\ \eta_i \end{bmatrix}
     \otimes \begin{bmatrix} 1 \\ Y_{i-1} \end{bmatrix} ,
 \end{align}
 provided that \ $n \sum_{i=1}^n Y_{i-1}^2 > \left(\sum_{i=1}^n Y_{i-1}\right)^2$.
\ By \eqref{help4}, the first factor converges almost surely to
 \[
   \bI_2
   \otimes
   \begin{bmatrix}
    1 & \EE(Y_\infty)\\
    \EE(Y_\infty) & \EE(Y_\infty^2)
   \end{bmatrix}^{-1}
   \qquad \text{as \ $n \to \infty$.}
 \]
For the second factor, we are going to apply the martingale central limit theorem (see
 Theorem \ref{Thm_CLT_discrete}) with the following choices: \ $d = 4$,
 \ $k_n = n$, \ $n \in \NN$, \ $\cF_{n,k} = \cF_k$, \ $n \in \NN$,
 \ $k \in \{1, \ldots, n\}$, \ and
 \[
   \bM_{n,k} = n^{-\frac{1}{2}}
            \sum_{i=1}^k
             \begin{bmatrix} \vare_i \\ \eta_i \end{bmatrix}
             \otimes
             \begin{bmatrix} 1 \\ Y_{i-1} \end{bmatrix} ,
   \qquad n \in \NN , \quad k \in \{1, \ldots, n\} .
 \]
Then, applying the identities \ $(\bA_1\otimes \bA_2)^\top = \bA_1^\top\otimes\bA_2^\top$ \ and
 \ $(\bA_1\otimes \bA_2)(\bA_3\otimes \bA_4)=(\bA_1\bA_3)\otimes(\bA_2\bA_4)$,
 \begin{align*}
  &\EE\big( (\bM_{n,k} - \bM_{n,k-1})(\bM_{n,k} - \bM_{n,k-1})^\top
            \,\big|\, \cF_{n,k-1} \big) \\
  &\qquad\qquad
   = \frac{1}{n}
     \EE\left( \left( \begin{bmatrix}
                       \vare_k \\ \eta_k
                      \end{bmatrix}
                      \otimes
                      \begin{bmatrix}
                       1 \\ Y_{k-1}
                      \end{bmatrix} \right)
               \left( \begin{bmatrix}
                       \vare_k \\ \eta_k
                      \end{bmatrix}
                      \otimes
                      \begin{bmatrix}
                       1 \\ Y_{k-1}
                      \end{bmatrix} \right)^\top
               \,\bigg|\, \cF_{k-1} \right) \\
  &\qquad\qquad
   = \frac{1}{n}
     \EE\left( \left( \begin{bmatrix}
                       \vare_k \\ \eta_k
                      \end{bmatrix}
                      \begin{bmatrix}
                       \vare_k \\ \eta_k
                      \end{bmatrix}^\top \right)
               \otimes
               \left( \begin{bmatrix}
                       1 \\ Y_{k-1}
                      \end{bmatrix}
                      \begin{bmatrix}
                       1 \\ Y_{k-1}
                      \end{bmatrix}^\top \right)
               \,\bigg|\, \cF_{k-1} \right) \\
  &\qquad\qquad
   = \frac{1}{n}
     \EE\left( \begin{bmatrix}
                \vare_k \\ \eta_k
               \end{bmatrix}
               \begin{bmatrix}
                \vare_k \\ \eta_k
               \end{bmatrix}^\top
               \,\bigg|\, \cF_{k-1} \right)
     \otimes
     \left( \begin{bmatrix}
             1 \\ Y_{k-1}
            \end{bmatrix}
            \begin{bmatrix}
             1 \\ Y_{k-1}
            \end{bmatrix}^\top \right) ,
  \qquad n\in\NN,\quad k\in\{1,\ldots,n\}.
 \end{align*}
Since \ $\EE(\vare_k^2 \mid \cF_{k-1}) = C_1 Y_{k-1} + C_2$, \ $k \in \NN$, \ and
 \ $\EE(\eta_k^2 \mid \cF_{k-1}) = C_3 Y_{k-1} + C_4$, \ $k \in \NN$, \ it
 remains to calculate
 \begin{align*}
  &\EE(\vare_k\eta_k \mid\cF_{k-1})
   = \EE\bigl( (Y_k - c - dY_{k-1})(X_k - X_{k-1} - \gamma - \delta Y_{k-1})
               \,\big|\, \cF_{k-1} \bigr) \\
  &= \EE\left( \left. \sigma_1
                      \int_{k-1}^k \ee^{-b(k-s)} \sqrt{Y_s} \, \dd W_s
                      \left( - \beta \sigma_1
                               \int_{k-1}^k \int_{k-1}^u
                                \ee^{-b(u-v)} \sqrt{Y_v} \, \dd W_v \, \dd u
                             + \sigma_2
                               \int_{k-1}^k \sqrt{Y_u} \, \dd \tW_u \right)
               \right| \cF_{k-1} \right) \\
  &= - \beta \sigma_1^2
       \int_{k-1}^k
        \EE\left( \left. \int_{k-1}^k \ee^{-b(k-s)} \sqrt{Y_s} \, \dd W_s
                         \int_{k-1}^u \ee^{-b(u-v)} \sqrt{Y_v} \, \dd W_v
                  \right| \cF_{k-1} \right) \dd u \\
  &\quad
     + \sigma_1 \sigma_2
       \EE\left( \left. \int_{k-1}^k \ee^{-b(k-s)} \sqrt{Y_s} \, \dd W_s
                        \int_{k-1}^k \sqrt{Y_u} \, \dd \tW_u
                 \right| \cF_{k-1} \right) .
\end{align*}
Again, by Equation (3.2.23) and Proposition 3.2.17 from Karatzas and Shreve
 \cite{KarShr}, we have
 \begin{align*}
  \EE(\vare_k \eta_k \mid \cF_{k-1})
  &= - \beta \sigma_1^2
       \int_{k-1}^k \int_{k-1}^u \ee^{-b(k+u-2v)}
        \EE(Y_v \mid \cF_{k-1}) \, \dd v \, \dd u \\
  &\quad
     + \sigma_1 \sigma_2 \varrho
       \int_{k-1}^k \ee^{-b(k-v)} \EE(Y_v \mid \cF_{k-1}) \, \dd v .
 \end{align*}
Using \eqref{cond_exp_discrete_Y}, by an easy calculation,
 \begin{align*}
  &\EE(\vare_k \eta_k \mid \cF_{k-1}) \\
  &= - \beta \sigma_1^2
       \int_{k-1}^k \int_{k-1}^u
        \ee^{-b(k+u-2v)}
        \left( \ee^{-b(v-k+1)} Y_{k-1}
               + a \int_{k-1}^v \ee^{-b(v-s)} \, \dd s \right)
        \dd v \, \dd u \\
  &\quad
     + \sigma_1 \sigma_2 \varrho
       \int_{k-1}^k
        \ee^{-b(k-v)}
        \left( \ee^{-b(v-k+1)} Y_{k-1}
               + a \int_{k-1}^v \ee^{-b(v-s)} \, \dd s \right)
        \dd v \\
  &= \left( - \beta \sigma_1^2
              \int_{0}^1 \int_{0}^{u'} \ee^{-b(u'-v'+1)} \, \dd v' \, \dd u'
            + \sigma_1 \sigma_2 \varrho \ee^{-b} \right)
     Y_{k-1}
     - a\beta \sigma_1^2
       \int_{0}^1 \int_{0}^{u'} \int_{0}^{v'}
        \ee^{-b(u'-v'-s'+1)} \, \dd s' \, \dd v' \, \dd u' \\
  &\quad
     + a\sigma_1 \sigma_2 \varrho
       \int_{0}^1 \int_{0}^{v'} \ee^{-b(1-s')} \, \dd s' \, \dd v'
   =: C_5 Y_{k-1} + C_6 , \qquad k \in \NN .
 \end{align*}
Hence, by \eqref{stac_moment} and \eqref{ergodic_discrete1},
 \begin{align*}
  &\sum_{k=1}^n
    \EE\big( (\bM_{n,k} - \bM_{n,k-1})(\bM_{n,k} - \bM_{n,k-1})^\top \mid \cF_{n,k-1} \big) \\
  &\qquad
   = \frac{1}{n}
     \sum_{k=1}^n
      \begin{bmatrix}
       C_1 Y_{k-1} + C_2 & C_5 Y_{k-1} + C_6 \\
       C_5 Y_{k-1} + C_6 & C_3 Y_{k-1} + C_4
      \end{bmatrix}
      \otimes
      \begin{bmatrix}
       1 & Y_{k-1} \\
       Y_{k-1} & Y_{k-1}^2
      \end{bmatrix} \\
  &\qquad
   = \frac{1}{n}
     \sum_{k=1}^n
      \begin{bmatrix}
       C_1 & C_5 \\
       C_5 & C_3
      \end{bmatrix}
      \otimes
      \begin{bmatrix}
       Y_{k-1} & Y_{k-1}^2 \\
       Y_{k-1}^2 & Y_{k-1}^3
      \end{bmatrix}
      + \frac{1}{n}
         \sum_{k=1}^n
          \begin{bmatrix}
           C_2 & C_6 \\
           C_6 & C_4
          \end{bmatrix}
          \otimes
          \begin{bmatrix}
           1 & Y_{k-1} \\
           Y_{k-1} & Y_{k-1}^2
          \end{bmatrix} \\
  &\qquad
   \as \begin{bmatrix}
        C_1 & C_5 \\
        C_5 & C_3
       \end{bmatrix}
       \otimes
       \begin{bmatrix}
        \EE(Y_\infty) & \EE(Y_\infty^2) \\
        \EE(Y_\infty^2) & \EE(Y_\infty^3)
       \end{bmatrix}
       + \begin{bmatrix}
          C_2 & C_6 \\
          C_6 & C_4
         \end{bmatrix}
         \otimes
         \begin{bmatrix}
          1 & \EE(Y_\infty) \\
          \EE(Y_\infty) & \EE(Y_\infty^2)
         \end{bmatrix}
   =: \bD
   \qquad \text{as \ $n \to \infty$,}
 \end{align*}
 where the $4\times 4$ limit matrix \ $\bD$ \ is necessarily
 symmetric and positive semi-definite (indeed, the limit of positive
 semi-definite matrices is positive semi-definite).

Next, we check Lindeberg condition \eqref{Lindeberg}.
Since
 \[
   \|\bx\|^2 \bbone_{\{\|\bx\|\geq\vare\}}
   \leq \frac{\|\bx\|^4}{\vare^2} \bbone_{\{\|\bx\|\geq\vare\}}
   \leq \frac{\|\bx\|^4}{\vare^2},
   \qquad \bx \in \RR^4 , \quad \vare \in \RR_{++} ,
 \]
 and \ $\|\bx\|^4=(x_1^2 + x_2^2 + x_3^2 + x_4^2)^2 \leq 4 (x_1^4 + x_2^4 + x_3^4 + x_4^4)$,
 \ $x_1, x_2, x_3, x_4 \in \RR$, \ it is enough to check that
 \begin{align*}
  &\frac{1}{n^2}
   \sum_{k=1}^n
    \bigl( \EE(\vare_k^4 \mid \cF_{k-1})
           + Y_{k-1}^4 \EE(\vare_k^4 \mid \cF_{k-1})
           + \EE(\eta_k^4 \mid \cF_{k-1})
           + Y_{k-1}^4 \EE(\eta_k^4 \mid \cF_{k-1}) \bigr) \\
  &\qquad\qquad
   = \frac{1}{n^2}
     \sum_{k=1}^n \EE( (1+Y_{k-1}^4)(\vare_k^4 + \eta_k^4) \mid \cF_{k-1})
   \stoch 0 \qquad \text{as \ $n \to \infty$.}
 \end{align*}
Instead of convergence in probability, we show convergence in \ $L^1$, \ i.e.,
 we check that
 \begin{align*}
  \frac{1}{n^2} \sum_{k=1}^n \EE((1+Y_{k-1}^4)(\vare_k^4 + \eta_k^4)) \to 0
  \qquad \text{as \ $n \to \infty$.}
 \end{align*}
Clearly, it is enough to show that
 \[
   \sup_{k\in\NN} \EE((1+Y_{k-1}^4)(\vare_k^4 + \eta_k^4)) < \infty .
 \]
By Cauchy--Schwarz inequality,
 \[
   \EE((1+Y_{k-1}^4)(\vare_k^4 + \eta_k^4))
   \leq \sqrt{\EE((1+Y_{k-1}^4)^2) \EE((\vare_k^4 + \eta_k^4)^2)}
   \leq \sqrt{2} \sqrt{\EE((1+Y_{k-1}^4)^2) \EE(\vare_k^8 + \eta_k^8)}
 \]
 for all \ $k \in \NN$.
\ Since, by Proposition 3 in Ben Alaya and Kebaier \cite{BenKeb2},
 \begin{align}\label{BenAlaya_Kebaier}
  \sup_{t\in\RR_+} \EE( Y_t^\kappa) < \infty, \qquad \kappa \in \RR_+ ,
 \end{align}
 it remains to check that \ $\sup_{k\in\NN} \EE(\vare_k^8 + \eta_k^8) < \infty$.
\ Since, by the power mean inequality,
 \[
   \EE(\vare_k^8)
   = \EE(\vert Y_k - d Y_{k-1} - c\vert^8)
   \leq \EE((Y_k + d Y_{k-1} + c)^8)
   \leq 3^7 \EE(Y_k^8 + d^8Y_{k-1}^8 + c^8) ,
   \qquad k \in \NN ,
 \]
 using \eqref{BenAlaya_Kebaier}, we have
 \ $\sup_{k\in\NN} \EE(\vare_k^8) < \infty$.
\ Using \eqref{Solutions} and again the power mean inequality, we have
 \begin{align*}
  \EE(\eta_k^8)
  &= \EE((X_k - X_{k-1} - \gamma - \delta Y_{k-1})^8) \\
  &= \EE\left( \left( \alpha - \beta \int_{k-1}^k Y_u \, \dd u
                      + \sigma_2 \varrho \int_{k-1}^k \sqrt{Y_u} \, \dd W_u
                      + \sigma_2 \sqrt{1-\varrho^2}
                        \int_{k-1}^k \sqrt{Y_u} \, \dd B_u
                      - \gamma - \delta Y_{k-1} \right)^8 \right) \\
  &\leq 6^7 \EE\Biggl( \alpha^8
                      + \beta^8 \left( \int_{k-1}^k Y_u \, \dd u \right)^8
                      + \sigma_2^8 \varrho^8
                        \left( \int_{k-1}^k \sqrt{Y_u} \, \dd W_u \right)^8
                      + \sigma_2^8 (1-\varrho^2)^4
                        \left( \int_{k-1}^k \sqrt{Y_u} \, \dd B_u \right)^8 \\
  &\phantom{ \leq 6^7 \EE\Biggl(}
                      + \delta^8 Y_{k-1}^8 + \gamma^8 \Biggr) ,
   \qquad k \in \NN .
 \end{align*}
By Jensen's inequality and \eqref{BenAlaya_Kebaier},
 \begin{align}\label{help8}
 \begin{split}
   \sup_{k\in\NN} \EE\left( \left( \int_{k-1}^k Y_u \, \dd u \right)^8 \right)
   &\leq \sup_{k\in\NN} \EE\left( \int_{k-1}^k Y_u^8 \, \dd u \right)
    = \sup_{k\in\NN} \int_{k-1}^k \EE(Y_u^8) \, \dd u\\
   &\leq \left(\sup_{t\in\RR_+} \EE(Y_t^8)\right) \left(\sup_{k\in\NN} \int_{k-1}^k 1\,\dd u\right)
         = \sup_{t\in\RR_+} \EE(Y_t^8)
   < \infty .
 \end{split}
 \end{align}
By the SDE \eqref{Heston_SDE} and the power mean inequality,
 \begin{align*}
  \EE\left(\left( \int_{k-1}^k \sqrt{Y_u} \, \dd W_u \right)^8\right)
   & \leq \frac{1}{\sigma_1^8} \EE\left( \left( Y_k - Y_{k-1} -a - b\int_{k-1}^k Y_u\,\dd u  \right)^8\right) \\
   &  \leq \frac{4^7}{\sigma_1^8} \EE\left( Y_k^8 + Y_{k-1}^8 + a^8 + b^8 \left( \int_{k-1}^k Y_u \, \dd u \right)^8 \right),
     \qquad k\in\NN,
 \end{align*}
 and hence, by \eqref{help8},
 \begin{align*}
   \sup_{k\in\NN} \EE\left( \left( \int_{k-1}^k \sqrt{Y_u} \, \dd W_u \right)^8 \right)
    \leq \frac{4^7}{\sigma_1^8} \left(2\sup_{t\in\RR_+} \EE(Y_t^8) + a^8 + b^8\sup_{t\in\RR_+} \EE(Y_t^8)\right)
     <\infty.
 \end{align*}
Further, using that the conditional distribution of \ $\int_{k-1}^k \sqrt{Y_u} \, \dd B_u$ \ given
 \ $(Y_u)_{u\in[0,k]}$ \ is normal with mean \ $0$ \ and variance \ $\int_{k-1}^k Y_u\,\dd u$ \ for all \ $k\in\NN$,
 \ we have
 \[
    \EE\left( \left(\int_{k-1}^k \sqrt{Y_u} \, \dd B_u \right)^8 \, \Big\vert \, (Y_u)_{u\in[0,k]} \right)
       = 105 \left(\int_{k-1}^k Y_u \, \dd u \right)^4,\qquad k\in\NN,
 \]
 and consequently
 \[
     \EE\left( \left(\int_{k-1}^k \sqrt{Y_u} \, \dd B_u \right)^8 \right)
        = 105 \EE\left( \left(\int_{k-1}^k Y_u \, \dd u \right)^4 \right),
        \qquad k\in\NN.
 \]
Hence, similarly to \eqref{help8}, we have
 \[
   \sup_{k\in\NN} \EE\left( \left(\int_{k-1}^k \sqrt{Y_u} \, \dd B_u \right)^8 \right)
      \leq 105 \sup_{t\in\RR_+} \EE(Y_t^4)
      <\infty,
 \]
 which yields that \ $\sup_{k\in\NN} \EE(\eta_k^8) < \infty$.
\ All in all, by the martingale central limit theorem (see, Theorem \ref{Thm_CLT_discrete}),
 \begin{align*}
  \bM_{n,n} =
  n^{-1/2}
  \sum_{k=1}^n
   \begin{bmatrix}
    \vare_k \\ \eta_k
   \end{bmatrix}
   \otimes
   \begin{bmatrix}
    1 \\ Y_{k-1}
   \end{bmatrix}
  \distr \cN_4\left( \bzero, \bD \right)
  \qquad \text{as \ $n \to \infty$.}
 \end{align*}

Consequently, by \eqref{help7} and Slutsky's lemma,
 \begin{align*}
  \sqrt{n}
  \begin{bmatrix}
   \hc_n^{\CLSED} - c \\
   \hd_n^{\CLSED} - d \\
   \hgamma_n^{\CLSED} - \gamma \\
   \hdelta_n^{\CLSED} - \delta
 \end{bmatrix}
 \distr
 \cN_4\left( \bzero,
             \left(\bI_2\otimes
             \begin{bmatrix}
              1 & \EE(Y_\infty) \\
              \EE(Y_\infty) & \EE(Y_\infty^2)
             \end{bmatrix}\right)^{-1} \bD
             \left(\bI_2\otimes
             \begin{bmatrix}
              1 & \EE(Y_\infty) \\
              \EE(Y_\infty) & \EE(Y_\infty^2)
             \end{bmatrix}\right)^{-1} \right)
 \end{align*}
 as \ $n \to \infty$, \ where the covariance matrix of the limit distribution
 takes the form
 \begin{align}
  &\left( \bI_2
          \otimes
          \begin{bmatrix}
           1 & \EE(Y_\infty) \\
           \EE(Y_\infty) & \EE(Y_\infty^2)
          \end{bmatrix}\right)^{-1}
   \bD
   \left( \bI_2
          \otimes
          \begin{bmatrix}
           1 & \EE(Y_\infty) \\
           \EE(Y_\infty) & \EE(Y_\infty^2)
          \end{bmatrix}\right)^{-1} \nonumber  \\
  &= \left(\begin{bmatrix}
           C_1 & C_5 \\
           C_5 & C_3 \\
         \end{bmatrix}
         \otimes
         \left(
         \begin{bmatrix}
             1 & \EE(Y_\infty) \\
             \EE(Y_\infty) & \EE(Y_\infty^2) \\
         \end{bmatrix}^{-1}
         \begin{bmatrix}
             \EE(Y_\infty) & \EE(Y_\infty^2) \\
             \EE(Y_\infty^2) & \EE(Y_\infty^3) \\
         \end{bmatrix}
         \right)
        \right)
    \left(\bI_2\otimes
           \begin{bmatrix}
             1 & \EE(Y_\infty) \\
             \EE(Y_\infty) & \EE(Y_\infty^2) \\
          \end{bmatrix}^{-1} \right) \nonumber  \\
  &\quad
    +\left(\begin{bmatrix}
           C_2 & C_6 \\
           C_6 & C_4 \\
         \end{bmatrix}
         \otimes
         \left(
         \begin{bmatrix}
             1 & \EE(Y_\infty) \\
             \EE(Y_\infty) & \EE(Y_\infty^2) \\
         \end{bmatrix}^{-1}
        \begin{bmatrix}
             1 & \EE(Y_\infty) \\
             \EE(Y_\infty) & \EE(Y_\infty^2) \\
         \end{bmatrix}
         \right)
        \right)
    \left(\bI_2\otimes
           \begin{bmatrix}
             1 & \EE(Y_\infty) \\
             \EE(Y_\infty) & \EE(Y_\infty^2) \\
          \end{bmatrix}^{-1} \right)\nonumber
  \end{align}
  \begin{align}
  & = \begin{bmatrix}
           C_1 & C_5 \\
           C_5 & C_3 \\
    \end{bmatrix}
    \otimes \left( \begin{bmatrix}
             1 & \EE(Y_\infty) \\
             \EE(Y_\infty) & \EE(Y_\infty^2) \\
         \end{bmatrix}^{-1}
         \begin{bmatrix}
             \EE(Y_\infty) & \EE(Y_\infty^2) \\
             \EE(Y_\infty^2) & \EE(Y_\infty^3) \\
         \end{bmatrix}
         \begin{bmatrix}
             1 & \EE(Y_\infty) \\
             \EE(Y_\infty) & \EE(Y_\infty^2) \\
         \end{bmatrix}^{-1}
         \right)  \nonumber  \\
  & \quad +
     \begin{bmatrix}
           C_2 & C_6 \\
           C_6 & C_4 \\
    \end{bmatrix}
    \otimes \left( \begin{bmatrix}
             1 & \EE(Y_\infty) \\
             \EE(Y_\infty) & \EE(Y_\infty^2) \\
         \end{bmatrix}^{-1}
         \begin{bmatrix}
             1 & \EE(Y_\infty) \\
             \EE(Y_\infty) & \EE(Y_\infty^2) \\
         \end{bmatrix}
         \begin{bmatrix}
             1 & \EE(Y_\infty) \\
             \EE(Y_\infty) & \EE(Y_\infty^2) \\
         \end{bmatrix}^{-1}
         \right)  \nonumber  \\
  & =\frac{1}{(\EE(Y_\infty^2) - (\EE(Y_\infty))^2)^2}
     \begin{bmatrix}
           C_1 & C_5 \\
           C_5 & C_3 \\
    \end{bmatrix}\nonumber \\
   & \phantom{=\;}
     \otimes
     \left( \begin{bmatrix}
             \EE(Y_\infty^2) & -\EE(Y_\infty) \\
             -\EE(Y_\infty) & 1 \\
           \end{bmatrix}
         \begin{bmatrix}
             \EE(Y_\infty) & \EE(Y_\infty^2) \\
             \EE(Y_\infty^2) & \EE(Y_\infty^3) \\
          \end{bmatrix}
         \begin{bmatrix}
             \EE(Y_\infty^2) & -\EE(Y_\infty) \\
             -\EE(Y_\infty) & 1 \\
         \end{bmatrix}
  \right) \nonumber  \\
 &\quad + \frac{1}{\EE(Y_\infty^2) - (\EE(Y_\infty))^2}
     \begin{bmatrix}
           C_2 & C_6 \\
           C_6 & C_4 \\
    \end{bmatrix}
    \otimes
  \begin{bmatrix}
    \EE(Y_\infty^2) & -\EE(Y_\infty) \\
    -\EE(Y_\infty) & 1 \\
  \end{bmatrix}  \nonumber\\
 &= \frac{1}{(\EE(Y_\infty^2) - (\EE(Y_\infty))^2)^2}
     \begin{bmatrix}
           C_1 & C_5 \\
           C_5 & C_3 \\
      \end{bmatrix} \nonumber  \\
 & \,\quad  \otimes
    \begin{bmatrix}
           -\EE(Y_\infty)( (\EE(Y_\infty^2))^2 - \EE(Y_\infty)\EE(Y_\infty^3) ) & (\EE(Y_\infty^2))^2 - \EE(Y_\infty)\EE(Y_\infty^3) \\
           (\EE(Y_\infty^2))^2 - \EE(Y_\infty)\EE(Y_\infty^3) & \EE(Y_\infty^3) - 2\EE(Y_\infty)\EE(Y_\infty^2) + (\EE(Y_\infty))^3 \\
         \end{bmatrix}  \nonumber \\
  &  \quad + \frac{1}{\EE(Y_\infty^2) - (\EE(Y_\infty))^2}
    \begin{bmatrix}
           C_2 & C_6 \\
           C_6 & C_4 \\
    \end{bmatrix}
    \otimes
  \begin{bmatrix}
    \EE(Y_\infty^2) & -\EE(Y_\infty) \\
    -\EE(Y_\infty) & 1 \\
  \end{bmatrix}  \nonumber
 \end{align}
 \begin{align}
 & =   \begin{bmatrix}
           C_1 & C_5 \\
           C_5 & C_3 \\
      \end{bmatrix}
      \otimes
      \begin{bmatrix}
           \frac{a(2a+\sigma_1^2)}{b\sigma_1^2} & -\frac{2a+\sigma_1^2}{\sigma_1^2} \\
           -\frac{2a+\sigma_1^2}{\sigma_1^2} & \frac{2b(a+\sigma_1^2)}{a\sigma_1^2} \\
         \end{bmatrix}
 +    \begin{bmatrix}
           C_2 & C_6 \\
           C_6 & C_4 \\
      \end{bmatrix}
       \otimes
     \begin{bmatrix}
           \frac{2a+\sigma_1^2}{\sigma_1^2} & -\frac{2b}{\sigma_1^2} \\
           -\frac{2b}{\sigma_1^2} & \frac{2b^2}{a\sigma_1^2} \\
         \end{bmatrix} := \bE.
   \label{matrixE}
 \end{align}
Indeed, by \eqref{stac_moment}, an easy calculation shows that
 \begin{align*}
  &\left(\EE(Y_\infty) \EE(Y_\infty^3) - (\EE(Y_\infty^2))^2\right) \EE(Y_\infty)
    = \frac{a^3\sigma_1^2}{4b^5}(2a+\sigma_1^2),\\
  &\EE(Y_\infty) \EE(Y_\infty^3) - (\EE(Y_\infty^2))^2
    = \frac{a^2\sigma_1^2}{4b^4}(2a+\sigma_1^2),\\
  &\EE(Y_\infty^3) - 2 \EE(Y_\infty) \EE(Y_\infty^2) + (\EE(Y_\infty))^3
    = \frac{a\sigma_1^2}{2b^3}(a+\sigma_1^2),\\
  &\EE(Y_\infty^2) - \left(\EE(Y_\infty)\right)^2
    = \frac{a\sigma_1^2}{2b^2}.
 \end{align*}
Finally, we show that \ $\bE$ \ is positive definite.
To show this, it is enough to check that
 \begin{itemize}
   \item[(i)] the matrix
         \[
          \begin{bmatrix}
            C_1 & C_5 \\
            C_5 & C_3 \\
          \end{bmatrix}
          \]
          is positive definite,
   \item[(ii)] the matrices
         \[
          \begin{bmatrix}
            C_2 & C_6 \\
            C_6 & C_4 \\
          \end{bmatrix},
          \qquad
          \begin{bmatrix}
           \frac{a(2a+\sigma_1^2)}{b\sigma_1^2} & -\frac{2a+\sigma_1^2}{\sigma_1^2} \\
           -\frac{2a+\sigma_1^2}{\sigma_1^2} & \frac{2b(a+\sigma_1^2)}{a\sigma_1^2} \\
         \end{bmatrix}
         \qquad \text{and} \qquad
         \begin{bmatrix}
           \frac{2a+\sigma_1^2}{\sigma_1^2} & -\frac{2b}{\sigma_1^2} \\
           -\frac{2b}{\sigma_1^2} & \frac{2b^2}{a\sigma_1^2} \\
         \end{bmatrix}
          \]
         are positive semi-definite.
 \end{itemize}
Indeed, the sum of a positive definite and a positive semi-definite square matrix is positive definite,
 the Kronecker product of positive semi-definite matrices is positive semi-definite and the Kronecker
 product of positive definite matrices is positive definite (as a consequence of the fact that the eigenvalues
 of the Kronecker product of two square matrices are the product of the eigenvalues of the two
 square matrices in question including multiplicities).
The positive semi-definiteness of the matrices
 \[
 \begin{bmatrix}
           \frac{a(2a+\sigma_1^2)}{b\sigma_1^2} & -\frac{2a+\sigma_1^2}{\sigma_1^2} \\
           -\frac{2a+\sigma_1^2}{\sigma_1^2} & \frac{2b(a+\sigma_1^2)}{a\sigma_1^2} \\
         \end{bmatrix}
 \qquad \text{and} \qquad
  \begin{bmatrix}
           \frac{2a+\sigma_1^2}{\sigma_1^2} & -\frac{2b}{\sigma_1^2} \\
           -\frac{2b}{\sigma_1^2} & \frac{2b^2}{a\sigma_1^2} \\
  \end{bmatrix}
 \]
 readily follows, since \ $\frac{a(2a+\sigma_1^2)}{b\sigma_1^2}>0$, \ $\frac{2a+\sigma_1^2}{\sigma_1^2}>0$, \ and the determinant of the
 matrices in question are \ $\frac{2a+\sigma_1^2}{\sigma_1^2}>0$ \ and \ $\frac{2b}{a\sigma_1^2}>0$, \ respectively.
Next, we prove that the matrices
  \[
  \begin{bmatrix}
            C_1 & C_5 \\
            C_5 & C_3 \\
   \end{bmatrix}
   \qquad \text{and} \qquad
   \begin{bmatrix}
            C_2 & C_4 \\
            C_4 & C_6 \\
   \end{bmatrix}
  \]
  are positive semi-definite.
Since \ $\PP(Y_0=y_0)=1$, \ we have \ $\EE(\vare_1^2 \mid \cF_0) = C_1 y_0 + C_2$,
 \ $\EE( \eta_1^2 \mid \cF_0) = C_3 y_0 + C_4$, \ and \ $\EE(\vare_1 \eta_1 \mid \cF_0) = C_5 y_0 + C_6$ \
 $\PP$-almost surely, hence
 \begin{align*}
  \EE( \vare_1^2 ) \EE(\eta_1^2)
   - \bigl( \EE( \vare_1 \eta_1 ) \bigr)^2
    = (C_1 C_3 - C_5^2) y_0^2 + (C_1 C_4 + C_2 C_3 - 2 C_5 C_6) y_0
      + C_2 C_4 - C_6^2  .
 \end{align*}
Clearly, by Cauchy--Schwarz's inequality,
 \[
  \EE( \vare_1^2 ) \EE( \eta_1^2 )
   - \bigl( \EE( \vare_1 \eta_1 ) \bigr)^2
   \geq 0 ,
 \]
 hence, by setting an arbitrary initial value \ $Y_0 = y_0 \in \RR_+$, \ we obtain
 \ $C_1 C_3 - C_5^2 \geq 0$ \ and \ $C_2 C_4 - C_6^2 \geq 0$.
\ Thus, both matrices
 \[
    \begin{bmatrix}
            C_1 & C_5 \\
            C_5 & C_3 \\
   \end{bmatrix}
   \qquad \text{and}\qquad
    \begin{bmatrix}
            C_2 & C_4 \\
            C_4 & C_6 \\
   \end{bmatrix}
  \]
 are positive semi-definite,
since \ $C_1 > 0$ \ and \ $C_2 > 0$.
\ Now we turn to check that
  \[
     \begin{bmatrix}
            C_1 & C_5 \\
            C_5 & C_3 \\
   \end{bmatrix}
  \]
  is positive definite.
Since \ $C_1 > 0$, \ this is equivalent to showing that $C_1 C_3 - C_5^2 > 0$.
\ Recalling the definition of the constants, we have
  \begin{align*}
 C_1 &= \sigma_1^2 \int_0^1 \ee^{-b(2-v)} \, \dd v = \sigma_1^2 \ee^{-2b}\frac{\ee^b - 1}{b}, \\
 C_3 &=
  \beta^2 \sigma_1^2
              \int_{0}^1 \int_{0}^1 \int_{0}^{u' \wedge v'}
               \ee^{-b(u'+v'-w')} \, \dd w' \, \dd v' \, \dd u'
              - 2 \beta \sigma_1 \sigma_2 \varrho
                \int_{0}^1 \int_{0}^{u'} \ee^{-bu'} \, \dd v' \, \dd u'
              + \sigma_2^2 \int_{0}^1 \ee^{-bu'} \, \dd u' \\
     &= b^{-3}\left(2 \ee^{-b} \beta^2 \sigma_1^2 (\sinh b -b) + 2 b \beta \varrho \sigma_1 \sigma_2((1 + b) \ee^{-b}-1) + b^2\sigma_2^2(1-\ee^{-b}) \right), \\
 C_5 &= - \beta \sigma_1^2
               \int_{0}^1 \int_{0}^{u'} \ee^{-b(u'-v'+1)} \, \dd v' \, \dd u'
             + \sigma_1 \sigma_2 \varrho \ee^{-b}
      = b^{-2} \sigma_1 \ee^{-b}\left(- \ee^{-b}\beta \sigma_1 (1+(b-1)\ee^{b}) + \varrho \sigma_2 b^2 \right),
  \end{align*}
  thus we have
  \begin{align*}
 C_1 C_3 - C_5^2 &= b^{-4}\ee^{-2b} \sigma_1^2
 	\Big(
 		2b(2+b^2)\beta \varrho \sigma_1 \sigma_2
  + 2 (\beta^2 \sigma_1^2 - 2 b \beta \varrho \sigma_1 \sigma_2 + b^2 \sigma_2^2) \cosh b	- (2+b^2)\beta^2 \sigma_1^2 \\
 		& \qquad \qquad \qquad - b^2 (2+b^2 \varrho^2) \sigma_2^2
 	\Big).
  \end{align*}
Consequently, using that $\cosh b = \sum_{k=0}^{\infty} \frac{b^{2k}}{(2k)!} > 1+\frac{b^2}{2}$ and that
 \[
 \beta^2 \sigma_1^2 - 2 b \beta \varrho \sigma_1 \sigma_2 + b^2 \sigma_2^2
    = (\beta \sigma_1 - b \varrho \sigma_2)^2 + b^2 (1-\varrho^2) \sigma_2^2 >0,
 \]
  we have
 \begin{align*}
 C_1 C_3 - C_5^2 &> b^{-4} \ee^{-2b} \sigma_1^2 \Big( 4b \beta \varrho \sigma_1 \sigma_2 + 2b^3 \beta \varrho \sigma_1 \sigma_2 + 2 \beta^2 \sigma_1^2 + b^2 \beta^2 \sigma_1^2 - 4 b \beta \varrho \sigma_1 \sigma_2 - 2 b^3 \beta \varrho \sigma_1 \sigma_2 \\
 & \qquad \qquad \qquad + 2 b^2 \sigma_2^2 + b^4 \sigma_2^2 - 2\beta^2 \sigma_1^2 - b^2 \beta^2 \sigma_1^2 - 2b^2 \sigma_2^2 - b^4 \varrho^2 \sigma_2^2\Big) \\
 &= b^{-4} \ee^{-2b} \sigma_1^2 ( b^4 (1-\varrho^2) \sigma_2^2 ) > 0.
 \end{align*}
 With this our proof is finished.
\proofend

So far we have obtained the limit distribution of the CLSE of the transformed parameters \ $(c,d,\gamma,\delta)$.
\ A natural estimator of \ $(a,b,\alpha,\beta)$ \ can be obtained from \eqref{CLSE_def}
 using relation \eqref{c_d_gamma_delta} detailed as follows.
Calculating the integrals in \eqref{c_d_gamma_delta} in the subcritical case, let us introduce
 the function \ $g: \RR_{++}^2 \times \RR^2 \to \RR_{++}\times (0,1)\times \RR^2$,
\begin{equation}\label{a_b_alpha_beta}
 g(a,b,\alpha,\beta):=
 \begin{bmatrix}
 a b^{-1} (1-\ee^{-b}) \\
 \ee^{-b} \\
 \alpha - a \beta b^{-2} (\ee^{-b}-1+b)\\
 -\beta b^{-1} (1-\ee^{-b}) \\
 \end{bmatrix}
 =\begin{bmatrix}
 c \\ d \\ \gamma \\ \delta
 \end{bmatrix},
 \qquad (a,b,\alpha,\beta)\in\RR_{++}^2 \times \RR^2.
 \end{equation}
Note that \ $g$ \ is bijective having inverse
 \begin{equation}\label{g_inverse}
  g^{-1}(c,d,\gamma,\delta)
  =
 \begin{bmatrix}
  - c \frac{\log d}{1-d} \\
 -\log d \\
 \gamma - c \delta \frac{d - 1 - \log d}{(1-d)^2} \\
 \delta \frac{\log d}{1-d}
 \end{bmatrix}
 =
 \begin{bmatrix}
  a \\ b \\ \alpha \\ \beta
 \end{bmatrix},
 \qquad (c,d,\gamma,\delta)\in\RR_{++}\times (0,1)\times \RR^2.
 \end{equation}
Indeed, for all \ $(c,d,\gamma,\delta)\in\RR_{++}\times (0,1)\times \RR^2$, \ we have
 \begin{align*}
    \alpha &= \gamma + a\beta b^{-2}(\ee^{-b} -1+b)
            = \gamma + (-c)\frac{\log d}{1-d}\delta \frac{\log d}{1-d} (-\log d)^{-2}(d -1-\log d)\\
           &= \gamma - c\delta \frac{d - 1 - \log d}{(1-d)^2}.
 \end{align*}
Under the conditions of Theorem \ref{c_d_normal} the CLSE
 \ $(\hc_n^{\CLSED}, \hd_n^{\CLSED}, \hgamma_n^{\CLSED}, \hdelta_n^{\CLSED})$ \ of \ $(c,d,\gamma,\delta)$ \
 is strongly consistent, hence in the subcritical case \ $(\hc_n^{\CLSED}, \hd_n^{\CLSED}, \hgamma_n^{\CLSED}, \hdelta_n^{\CLSED})$ \
 fall into the set \ $\RR_{++}\times (0,1)\times \RR^2$ \ for sufficiently large \ $n\in\NN$ \ with probability one.
Hence, in the subcritical case, one can introduce a natural estimator of
 \ $(a,b,\alpha,\beta)$ \ based on discrete time observations
 \ $(Y_i, X_i)_{i\in\{1,\ldots,n\}}$ \ by applying the inverse of \ $g$ \ to the
 CLSE of \ $(c,d,\gamma,\delta)$,\
 i.e.,
 \begin{equation}\label{CLSEabalphabeta_discrete}
  (\ha_n, \hb_n, \halpha_n, \hbeta_n)
  := g^{-1}(\hc_n^{\CLSED},\hd_n^{\CLSED},\hgamma_n^{\CLSED},\hdelta_n^{\CLSED})
 \end{equation}
 for sufficiently large \ $n\in\NN$ \ with probability one.
\begin{Rem}
We would like to stress the point that the estimator of
 \ $(a,b,\alpha,\beta)$ \ introduced in \eqref{CLSEabalphabeta_discrete}
 exists only for sufficiently
 large \ $n\in\NN$ \ with probability of \ $1$.
\ However, as all our results are asymptotic, this will not cause a problem.
From the considerations before this remark, we obtain
 \begin{equation}\label{CLSE_abab}
   \bigl(\ha_n, \hb_n, \halpha_n, \hbeta_n\bigr)
   = \argmin_{(a,b,\alpha,\beta)\in\RR_{++}^2\times\RR^2}
      \sum_{i=1}^n
       \left[ (Y_i - d Y_{i-1} - c)^2
              + (X_i - X_{i-1} - \gamma - \delta Y_{i-1})^2 \right]
 \end{equation}
 for sufficiently large \ $n\in\NN$ \ with probability one.
We call the attention that \ $\bigl(\ha_n, \hb_n, \halpha_n, \hbeta_n\bigr)$
 \ does not necessarily provides a CLSE of \ $(a,b,\alpha,\beta)$, \ since
 in \eqref{CLSE_abab} one takes the infimum only on the set
 \ $\RR_{++}^2\times\RR^2$ \ instead of \ $\RR^4$.
\ Formula \eqref{CLSE_abab} serves as a motivation for calling
 \ $\bigl(\ha_n, \hb_n, \halpha_n, \hbeta_n\bigr)$ \ essentially conditional least
 squares estimator in the Abstract.
\proofend
\end{Rem}

\begin{Thm}\label{Thm_CLSE_a_b_alpha_beta}
Under the conditions of Theorem \ref{c_d_normal} the sequence
 \ $\bigl(\ha_n, \hb_n, \halpha_n, \hbeta_n\bigr)$, \ $n \in \NN$, \ is
 strongly consistent and asymptotically normal, i.e.,
\[
   (\ha_n, \hb_n, \halpha_n, \hbeta_n)
   \as (a, b, \alpha, \beta) \qquad \text{as \ $n \to \infty$,}
 \]
 and
\[
\sqrt{n}
\begin{bmatrix}
\ha_n - a \\
\hb_n - b \\
\halpha_n - \alpha \\
\hbeta_n - \beta
\end{bmatrix}
\distr
\cN_4
\left(
0,\bJ \bE \bJ^\top \right)
      \qquad \text{as \ $n\to\infty$,}
\]
 where  \ $\bE \in \RR^{2\times 2}$ \ is a symmetric, positive definite matrix given in \eqref{matrixE} and
\[
\bJ :=
\begin{bmatrix}
-\frac{\log d}{1-d} & -c \frac{\log d - 1 + d^{-1}}{(1-d)^2} & 0 & 0 \\
0 & -\frac{1}{d} & 0 & 0 \\
 \delta \frac{\log d + 1 - d}{(1-d)^2}  &
 c \delta \frac{2\log d - d + d^{-1} }{(1-d)^3} & 1 & c \frac{\log d + 1 - d}{(1-d)^2} \\
0 & \delta \frac{\log d - 1 + d^{-1}}{(1-d)^2} & 0 & \frac{\log d}{1-d}
\end{bmatrix}
\]
with \ $c$, \ $d$ \ and \ $\delta$ \ given in
 \eqref{c_d_gamma_delta}.
\end{Thm}

{\noindent \bf Proof.}
The strong consistency of \ $(\ha_n, \hb_n, \halpha_n, \hbeta_n)$,
 \ $n \in \NN$, \ follows from the strong consistency of the CLSE of \ $(c,d,\gamma,\delta)$ \
  proved in Theorem \ref{c_d_normal} using also that the inverse function \ $g^{-1}$ \ given in
 \eqref{g_inverse} is continuous on \ $\RR_{++}\times (0,1)\times \RR^2$.
\ For the second part of the theorem we use Theorem \ref{c_d_normal}, and the so-called delta method
 (see, e.g., Theorem 11.2.14 in Lehmann and Romano \cite{LehRom}).
Indeed, one can extend the function \ $g^{-1}$ \ to be defined on \ $\RR^4$ \ not only on
 \ $\RR_{++}\times (0,1)\times \RR^2$ \ (e.g., let it be zero on the complement of
 \ $\RR_{++}\times (0,1)\times \RR^2$), \ $(\ha_n, \hb_n, \halpha_n, \hbeta_n)$ \ takes the form given
 in \eqref{CLSEabalphabeta_discrete} with this extension of \ $g^{-1}$ \ as well, and the Jacobian of
 \ $g^{-1}$ \ at \ $(c,d,\gamma,\delta)\in \RR_{++}\times (0,1)\times \RR^2$ \ is clearly \ $\bJ$.
\proofend

\section*{Acknowledgements}
We are undoubtedly grateful for the referee for pointing out mistakes in the proof of Theorem \ref{c_d_normal},
 and also for his/her several valuable comments that have led to an improvement of the manuscript.


\begin{thebibliography}{99}

\bibitem{BarPap2}
\textsc{Barczy, M.} and \textsc{Pap, G.} (2016).
Asymptotic properties of maximum-likelihood estimators for Heston models based on continuous time observations.
\textit{Statistics}
\textbf{50(2)} 389--417.

\bibitem{BarDorLiPap2}
\textsc{Barczy, M.}, \textsc{D\"oring, L.}, \textsc{Li, Z.} and
 \textsc{Pap, G.} (2014).
Stationarity and ergodicity for an affine two-factor model.
\textit{Advances in Applied Probability}
\textbf{46(3)} 878--898.

\bibitem{BenKeb2}
\textsc{Ben Alaya, M.} and \textsc{Kebaier, M.} (2013).
Asymptotic behavior of the maximum likelihood estimator for ergodic and
 nonergodic square-root diffusions.
\textit{Stochastic Analysis and Applications}
\textbf{31(4)} 552--573.

\bibitem{CoxIngRos}
\textsc{Cox, J. C.}, \textsc{Ingersoll, J. E.} and \textsc{Ross, S. A.} (1985).
A theory of the term structure of interest rates.
\textit{Econometrica}
\textbf{53(2)} 385--407.

\bibitem{Dok}
\textsc{Dokuchaev, N.} (2015).
 Pathwise estimation of the diffusion term for Cox-Ingersoll-Ross and similar processes.
 Available on the ArXiv: \texttt{http://arxiv.org/abs/1506.05627}

\bibitem{Dud}
\textsc{Dudley, R. M.} (1989).
\textit{Real Analysis and Probability}.
Wadsworth \& Brooks/Cole Advanced Books \& Software, Pacific Grove, California.

\bibitem{Hes}
\textsc{Heston, S.} (1993).
A closed-form solution for options with stochastic volatilities with
 applications to bond and currency options.
\textit{The Review of Financial Studies}
\textbf{6} 327--343.

\bibitem{HurLinMcC}
\textsc{Hurn, A. S., Lindsay, K. A.} and \textsc{McClelland, A. J.} (2013).
A quasi-maximum likelihood method for estimating the parameters of
 multivariate diffusions.
\textit{Journal of Econometrics}
\textbf{172} 106--126.

\bibitem{JSh}
\textsc{Jacod, J.} and \textsc{Shiryaev, A. N.} (2003).
\textit{Limit Theorems for Stochastic Processes}, 2nd ed.
Springer-Verlag, Berlin.

\bibitem{JinManRudTra}
\textsc{Jin, P., Mandrekar, V., R\"udiger, B.} and \textsc{Trabelsi, C.} (2013).
Positive Harris recurrence of the CIR process and its applications.
\textit{Communications on Stochastic Analysis}
\textbf{7(3)} 409--424.

\bibitem{JinRudTra}
\textsc{Jin, P., R\"udiger, B.} and \textsc{Trabelsi, C.}  (2016).
Positive Harris recurrence and exponential ergodicity of the basic affine jump-diffusion.
\textit{Stochastic Analysis and Applications}
\textbf{34(1)} 75--95.

\bibitem{KarShr}
\textsc{Karatzas, I.} and \textsc{Shreve, S. E.} (1991).
\textit{Brownian Motion and Stochastic Calculus,} 2nd ed.
 Springer-Verlag.

\bibitem{LehRom}
\textsc{Lehmann, E.L.} and \textsc{Romano, J. P.} (2009).
\textit{Testing Statistical Hypotheses,} 3rd ed.
 Springer-Verlag Berlin Heidelberg.

\bibitem{LiMa}
\textsc{Li, Z.} and \textsc{Ma, C.} (2015).
 Asymptotic properties of estimators in a stable Cox-Ingersoll-Ross model.
 \textit{Stochastic Processes and their Applications}
 \textbf{125(8)} 3196--3233.


\bibitem{OveRyd}
\textsc{Overbeck, L.} and \textsc{Ryd\'en, T.} (1997).
Estimation in the Cox-Ingersoll-Ross model.
\textit{Econometric Theory}
\textbf{13(3)} 430--461.

\bibitem{Shi}
\textsc{Shiryaev, A. N.} (1989).
\textit{Probability}, 2nd ed. Springer, New York.

\end{thebibliography}
\end{document}